\documentclass[12pt, reqno]{amsart}
\usepackage{ifthen,amsfonts,amsmath,amssymb,epic,eepic}
\usepackage{epsfig,color}

\makeatother

\theoremstyle{definition}

\def\fnum{equation}
\newtheorem{Thm}[\fnum]{Theorem}
\newtheorem{Cor}[\fnum]{Corollary}

\newtheorem{Lem}[\fnum]{Lemma}

\newtheorem{Pro}[\fnum]{Proposition}

\newcommand{\Vol}{{\text{Vol}}}

\newcommand{\nn}{{\bf{n}}}

\newcommand{\dist}{{\text {dist}}}

\def\ZZ{{\bold Z}}
\def\RR{{\bold R}}

\def\CC{{\bold C }}

\newcommand{\e}{{\text {e}}}
\newcommand{\Area}{{\text {Area}}}

\newcommand{\cF}{{\mathcal{F}}}

\newcommand{\cD}{{\mathcal{D}}}
\newcommand{\K}{{\text{K}}}

\newcommand{\cB}{{\mathcal{B}}}

\newcommand{\cP}{{\mathcal{P}}}
\newcommand{\cS}{{\mathcal{S}}}

\newcommand{\cA}{{\mathcal{A}}}

\newcommand{\eqr}[1]{(\ref{#1})}

\newcommand{\vb}{V \negmedspace B}

\newcommand{\cone}{{\bf{C}}}

\begin{document}

\title[Embedded minimal disks]
{Embedded minimal disks}

\author{Tobias H. Colding}%
\address{Courant Institute of Mathematical Sciences and Princeton University\\
251 Mercer Street\\ New York, NY 10012 and Fine Hall, Washington
Rd., Princeton, NJ 08544-1000}
\author{William P. Minicozzi II}%
\address{Department of Mathematics\\
Johns Hopkins University\\
3400 N. Charles St.\\
Baltimore, MD 21218}
\thanks{The first author was partially supported by NSF Grant DMS 9803253
and an Alfred P. Sloan Research Fellowship and the second author
by NSF Grant DMS 9803144 and an Alfred P. Sloan Research
Fellowship.}


\email{colding@cims.nyu.edu, minicozz@jhu.edu}


\maketitle

{\small
\tableofcontents}

\numberwithin{equation}{section}


\setcounter{part}{0}
\numberwithin{section}{part} 
\renewcommand{\rm}{\normalshape} 
\renewcommand{\thepart}{\Roman{part}}
\setcounter{section}{1}

\part{The main theorem - the limit foliation and the singular curve}
\label{p:p1}

This paper is a survey of our results about embedded minimal
disks.  Unlike our expository article \cite{CM15}, this paper is
intended for readers with some background knowledge on minimal
surfaces. However, even an expert reader may find it worthwhile to
look at \cite{CM15} first.

Let us start with two key examples of embedded minimal disks in
$\RR^3$.  The first are minimal graphs over  simply connected
domains and the second is the double spiral staircase (cf.
\cite{CM15}) known as the helicoid.

\vskip2mm \noindent {\bf{Example 1}}: Graphs of functions
$u:\Omega\to \RR$ where $\Omega\subset \RR^2$ is simply connected
and $u$ satisfies the second order nonlinear elliptic equation in
divergence form (the so-called   minimal surface equation)
\begin{equation}  \label{e:mineq}
\text{div} \left( \frac{\nabla u}{\sqrt{1+|\nabla u|^2}}\right)=0\, .
\end{equation}

\vskip2mm
\noindent
{\bf{Example 2}}:
(Helicoid).  See fig. \ref{f:f1}.  The minimal
surface in $\RR^3$ parametrized by
\begin{equation}
(s\cos t,s\sin t,-t)\text{ where }s,\,t\in \RR\, .
\end{equation}

\begin{figure}[htbp]
    \setlength{\captionindent}{20pt}
    \begin{minipage}[t]{0.5\textwidth}
    \centering\input{unot1.pstex_t}
    \caption{The helicoid is obtained
    by gluing together two $\infty$-valued graphs along a line.
The two multi-valued graphs are given in polar coordinates
by $u_1(\rho,\theta)=-\theta$ and $u_2(\rho,\theta)=-\theta+\pi$.
In either case $w(\rho,\theta)=-2\,\pi$.}\label{f:f1}
    \end{minipage}\begin{minipage}[t]{0.5\textwidth}
    \centering\input{unot2a.pstex_t}
    \caption{The separation of a multi-valued graph.  (Here the
multi-valued graph is shown with negative separation.)}\label{f:f2}
    \end{minipage}%
\end{figure}

\vskip2mm One of our main theorems is that every embedded minimal
disk can either be modelled by a minimal graph or by a piece of
the helicoid depending on whether the curvature is small or not;
see Theorem \ref{t:t0.1} below.  We will in this survey discuss
some of the ingredients in the proof of this and how those
ingredients fit together.

To be able to discuss the helicoid some more and in particular
give a precise meaning to that it is like a double spiral
staircase, we will need the notion of a multi-valued graph; see
fig. \ref{f:f2}. Let $D_r$ be the disk in the plane centered at
the origin and of radius $r$ and let $\cP$ be the universal cover
of the punctured plane $\CC\setminus \{0\}$ with global polar
coordinates $(\rho, \theta)$ so $\rho>0$ and $\theta\in \RR$.  An
$N$-valued graph of a function $u$ on the annulus $D_s\setminus
D_r$ is a single valued graph over
\begin{equation}
    \{(\rho,\theta)\,|\,r\leq \rho\leq s\, ,\, |\theta|\leq
N\,\pi\} \, .
\end{equation}
  The multi-valued graphs that we will
consider will never close up; in fact they will all be embedded.
Note that embedded corresponds to that the separation never
vanishes.  Here the separation (see fig. \ref{f:f2}) is the
function given by
\begin{equation}
w(\rho,\theta)=u(\rho,\theta+2\pi)-u(\rho,\theta)\, .
\end{equation}
If $\Sigma$ is the helicoid, then
$\Sigma\setminus x_3-\text{axis}=\Sigma_1\cup \Sigma_2$,
 where $\Sigma_1$, $\Sigma_2$ are $\infty$-valued graphs.
$\Sigma_1$ is the graph of the function
$u_1(\rho,\theta)=-\theta$ and $\Sigma_2$ is the
graph of the function $u_2(\rho,\theta)=-\theta+\pi$.
In either case the separation
$w=-2\,\pi$.  A multi-valued minimal graph is a multi-valued graph of a
function $u$ satisfying the minimal surface equation \eqr{e:mineq}.

Here, as in \cite{CM6} and \cite{CM8}, we have normalized so our embedded
multi-valued graphs have negative separation.  This can be achieved
after possibly reflecting in a plane.

\vskip6mm
Let now $\Sigma_i\subset B_{2R}$ be a
sequence of embedded minimal
 disks with
$\partial \Sigma_i\subset \partial B_{2R}$.  Clearly
(after possibly going to a subsequence) either (1) or (2) occur:\\
(1) $\sup_{B_{R}\cap\Sigma_i}|A|^2\leq C<\infty$ for some constant $C$.\\
(2) $\sup_{B_{R}\cap\Sigma_i}|A|^2\to \infty$.\\
In (1) (by a standard argument) the intrinsic ball $\cB_s(y_i)$ is a
graph for all $y_i\in B_{R}\cap \Sigma_i$, where $s$ depends only on $C$.
Thus the main case is (2) which is the subject of the next theorem.

\vskip6mm
Using the notion of multi-valued graphs, this our main
theorem, can now be stated:

\begin{figure}[htbp]
    \setlength{\captionindent}{20pt}
    \begin{minipage}[t]{0.5\textwidth}
    \centering\input{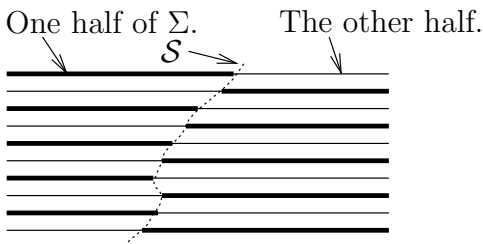}
    \caption{Theorem \ref{t:t0.1} - the singular set, $\cS$, and
the two multi-valued graphs.}\label{f:f3}
    \end{minipage}
\end{figure}

\begin{Thm} \label{t:t0.1}
(Theorem 0.1 in \cite{CM6}).  See fig. \ref{f:f3}.
Let $\Sigma_i \subset B_{R_i}=B_{R_i}(0)\subset \RR^3$ be a
sequence of embedded minimal
disks with $\partial \Sigma_i\subset \partial B_{R_i}$
where $R_i\to \infty$. If $\sup_{B_1\cap \Sigma_i}|A|^2\to \infty$, then
there exists a subsequence, $\Sigma_j$,
and
a Lipschitz curve $\cS:\RR\to \RR^3$ such that after a rotation of $\RR^3$:\\
\underline{1.} $x_3(\cS(t))=t$.  (That is, $\cS$ is a graph
over the $x_3$-axis.)\\
\underline{2.}  Each $\Sigma_j$ consists of exactly two multi-valued
graphs away
from $\cS$ (which spiral together).\\
\underline{3.} For each $1>\alpha>0$, $\Sigma_j\setminus \cS$ converges
in the $C^{\alpha}$-topology to the foliation,
$\cF=\{x_3=t\}_t$, of $\RR^3$.\\
\underline{4.}  $\sup_{B_{r}(\cS (t))\cap \Sigma_j}|A|^2\to\infty$ for
all $r>0$, $t\in \RR$.  (The curvatures blow up along $\cS$.)
\end{Thm}

In \underline{2.}, \underline{3.}
that $\Sigma_j\setminus \cS$ are multi-valued graphs and
converges to $\cF$ means that for each compact subset
$K\subset \RR^3\setminus \cS$
and $j$ sufficiently large $K\cap \Sigma_j$ consists of multi-valued
graphs over (part of) $\{x_3=0\}$
and $K\cap \Sigma_j\to K\cap \cF$ in the sense of graphs.

Theorem \ref{t:t0.1} (as many of the other results discussed
below) is modelled by the helicoid and its rescalings. Take a
sequence $\Sigma_i = a_i \, \Sigma$ of rescaled helicoids where
$a_i \to 0$. The curvatures of this sequence are blowing up along
the vertical axis. The sequence converges (away from the vertical
axis) to a foliation by flat parallel planes. The singular set
$\cS$ (the axis) then consists of removable singularities.

\vskip2mm
Of the many different paths that one could choose through our results about
embedded minimal disks, we have chosen here one which discusses some of our
more elementary results in greater detail and then only gives a very rough
overview of some of our key results which are more difficult.  Our hope is
that by doing so this survey can serve as a reading guide
for our papers where the
reader is eased into the subject and then is shown the anatomy of the
proof of our main theorem.

One of these more elementary themes is that of analysis of
multi--valued embedded minimal graphs.  From estimates on the
growth and decay of the separation between consecutive sheets of
such graphs to the results about why such graphs are proper in a
certain qualitative sense if they are contained in larger embedded
minimal disks.  These results about multi--valued graphs are among
what is discussed in some detail in the first part where we also
discuss some of the other main results that go into the proof of
our main theorem that any embedded minimal disk is either a graph
of a function or can be approximated by a piece of a helicoid. In
addition, we also discuss in some detail in the first part why the
singular set, that is the set of points of large curvature in an
embedded minimal disk, must all be contained in a curve which is a
Lipschitz graph over a straight line.

In the second part we come to a  less elementary result. Namely,
that of why near a point of large curvature of an embedded minimal
disk there must inside the disk be a small multi--valued graph.
Here small means on the scale of one over the square--root of the
maximum of the curvature.   The much less elementary, but key,
result of why such small multi--valued graphs extend to large ones
is only discussed very briefly in the second section of that
second part.

In the third part  we discuss the crucial one--sided curvature
estimate. This is the estimate that gives a cone condition for all
the points of large curvature.  Namely, given a point of large
curvature of an embedded minimal disk, then the cone condition is
that all the other points of large curvature must lie within a
double convex cone of the initial point of large curvature.
Iterating this condition gives that the set of all the points of
large curvature is the Lipschitz graph over a straight line.  The
one--sided curvature estimates uses all the results discussed in
the first part.

In the final part we discuss what the differences are between the
so-called local and global case.  The local case is where we have
a sequence of embedded minimal disks in a ball of fixed radius in
$\RR^3$ - the global case is where the disks are in a sequence of
expanding balls with radii tending to infinity.  As we will see in
the final part, then in the local case we can get limits with
singularities. In the global case this does not happen because in
fact any limit is always a foliation by flat parallel planes (this
is assuming that the curvatures of the sequence are blowing up).

There are a number of important results that go into the proof of
our main theorem about embedded minimal disks that are either not
discussed here or are barely mentioned.  One of these is why given
a point of large curvature in an embedded minimal disk there are
points of large curvature nearby above and below.  This was one of
the key results proven in \cite{CM5} but is quite technical and
thus has been omitted from this survey.

\vskip6mm  Let $x_1 , x_2 , x_3$ be the standard coordinates on
$\RR^3$. For $y \in \Sigma \subset \RR^3$ and $s > 0$, the
extrinsic and intrinsic balls are $B_s(y)$, $\cB_s(y)$.
$\K_{\Sigma}$ the sectional curvature of a smooth compact surface
$\Sigma$ and when $\Sigma$ is immersed $A_{\Sigma}$ will be its
second fundamental form (so when $\Sigma$ is minimal, then
$|A|^2=-2\,\K_{\Sigma}$). When $\Sigma$ is oriented,
$\nn_{\Sigma}$ is the unit normal.

\vskip6mm
Using Theorems \ref{t:t0.1}, \ref{t:t2},
W. Meeks and H. Rosenberg proved
that the plane and helicoid are the only
complete properly embedded
simply-connected minimal surfaces in $\RR^3$, \cite{MeRo}.

\section{Estimates for multi-valued minimal graphs} \label{s:s0}

We will later see that, just like the helicoid, general embedded
minimal disks with large curvature at some interior point can be
built out of multi-valued graphs. This will be particularly useful
once we have a good understanding of general embedded multi-valued
minimal graphs. In this section, we discuss three basic, but
useful, estimates for such graphs:
\begin{itemize}
\item
A gradient estimate for the separation which implies sublinear
growth of $w$ when there are enough sheets.
\item
A curvature estimate for $2$-valued minimal graphs whose
separation grows sublinearly.
\item
Sharp logarithmic upper and lower bounds for the separation when
there is a growing number of sheets.
\end{itemize}
(See also \cite{CM7} for further discussion of these estimates,
their analogs for minimal annuli, and their
implications.)

\vskip2mm The first of these basic estimates was obtained in
\cite{CM3} where we showed (essentially by a gradient estimate)
that the separation between the sheets grows sublinearly (in
Theorem \ref{t:main2} below we will discuss an improvement of this
when the number of sheets grows sufficiently fast). Precisely, in
\cite{CM3} we showed that:

\begin{Pro}
(Proposition II.2.12 in \cite{CM3}). Given $\alpha > 0$, there
exists $N_{\alpha}$  so if $\Sigma$ is an embedded
$N_{\alpha}$-valued minimal graph over $D_{\e^{N_{\alpha}} \, R}
\setminus D_{\e^{-N_{\alpha}}}$ of $u$ and $1\leq \rho \leq R$,
then
\begin{equation}    \label{e:slg}
    \rho^{-\alpha} \leq \frac{w(\rho,0)}{ w(1,0)}\leq \rho^{\alpha}  \, .
\end{equation}
\end{Pro}

Thus, choosing $\alpha<1$, \eqr{e:slg} gives the sublinear growth
of the separation $w$.
 This sublinear growth  of the separation  is the main
 benefit of having at least $N_{\alpha}$ sheets.
 It comes from integrating the gradient bound proven in Proposition II.2.12 in \cite{CM3}
\begin{equation}   \label{e:wantit}
   \frac{ |\nabla  w |}{|w|} \leq  \frac{\alpha}{\rho}  \, .
\end{equation}
 Many of our
 results on embedded multi-valued graphs apply as long as \eqr{e:slg} holds
 and we have at least two sheets.

\vskip2mm To get the better logarithmic bounds on the separation,
one needs curvature estimates for embedded multi-valued graphs. It
follows from Heinz's curvature estimate for minimal graphs
(theorem $2.4$ in \cite{CM1}) that (away from its boundary) a
multi-valued minimal graph has quadratic curvature decay
\begin{equation}  \label{e:qcd}
|A|^2 \leq C \, r^{-2} \, .
\end{equation}
This scale-invariant estimate \eqr{e:qcd}  does not require
embeddedness  and can easily be seen to be sharp without any
further assumptions on the graph. However, using the embeddedness
-- and, in particular, the sublinear growth of the separation that
embeddedness implies --
 we showed in
corollary 2.3 of \cite{CM8} that the curvature of a multi-valued
embedded minimal graph decays faster than quadratically.  That is,
we showed that
\begin{equation}  \label{e:5/6}
|A|^2 \leq C \, r^{-2-5/6}
\end{equation}
for embedded $2$-valued minimal graphs whose separation grows
sublinearly; cf. \eqr{e:slg}.  (Note that by \cite{CM3}, there
exists $N_{\alpha}$ so this applies to any
 $N_{\alpha}$-valued embedded minimal graph.)  The important point for the
applications in \cite{CM8} is that the faster than linear decay on
the Hessian of $u$ that \eqr{e:5/6} implies together with the
nonlinear form of \eqr{e:mineq} gives that $|\Delta\,u|$ decays
faster than quadratically; see equation (3.6) in \cite{CM8}.

The above faster than quadratic bound \eqr{e:5/6}
used the nonlinearity of the minimal
surface equation much like the nonlinearity is needed for proving
the Bernstein theorem.  In other situations the nonlinearity seems
more to add difficulties than being of any help. In those
situations one tries to model the minimal surface equation by the
Laplace equation and use that from \eqr{e:5/6} if $u$ is a
multi-valued solution of the minimal surface equation, then
$|\Delta\,u|$ decays faster than quadratically so that $u$
``become more and more like a harmonic function''. We will now
take advantage of this to discuss some analysis of such
multi-valued solutions that we will need later.  This analysis is
from \cite{CM8}.

The first such result is the following sharp upper and lower
logarithmic bound
for the separation of a multi-valued graph:

\begin{figure}[htbp]
    \setlength{\captionindent}{20pt}
    \begin{minipage}[t]{0.5\textwidth}
    \centering\input{dis4.pstex_t}
    \caption{Theorem \ref{t:main2}:
The sharp logarithmic upper and lower bounds for the separation.}
\label{f:f4}
\end{minipage}\begin{minipage}[t]{0.5\textwidth}
\centering\input{dis5.pstex_t}
    \caption{The domain of $u$ in Theorem \ref{t:main2}
in $(\log \rho,\theta)$-coordinates.}\label{f:f5}
    \end{minipage}
\end{figure}

\begin{Thm}  \label{t:main2}
(Theorem 1.3 of \cite{CM8}).
See fig. \ref{f:f4} and fig. \ref{f:f5}.
Given $c>0$, there exists $c_1$, $c_2$ such that
if $u$ satisfies \eqr{e:mineq} on
$\{(\rho,\theta)\,|1/2\leq \rho\leq R\text{ and }
\,c\,|\theta|\leq (1 + \log \rho )^{3/2} \}$ and $w<0$
together with a  slight condition on the growth of $u$ and $w$
(see equation (1.5) in \cite{CM8}),
then  for $2\leq \rho \leq R^{1/2}$
\begin{equation}    \label{e:upandlow}
\frac{c_1}{\log \rho}\leq \frac{w(\rho,0)}{ w(1,0)}\leq c_2 \log \rho\, .
\end{equation}
\end{Thm}

\begin{figure}[htbp]
    \setlength{\captionindent}{20pt}
    \begin{minipage}[t]{0.5\textwidth}
    \centering\input{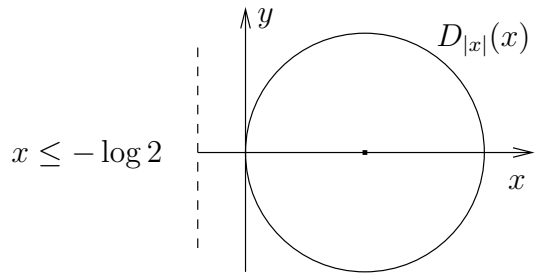}
    \caption{Proof of Theorem \ref{t:main2}:  Applying the mean
value equality to a negative harmonic function defined in a half-space.}
\label{f:f6}
    \end{minipage}
\end{figure}

{\bf{The idea of the lower bound for the separation}}, i.e.,
$w(\rho,0)/w(1,0)\geq c_1/\log\rho$:  If $u$ is as in Theorem
\ref{t:main2}, then $w<0$ is almost harmonic and it can be seen to
follow from \eqr{e:5/6} and the form of \eqr{e:mineq} that so is
(the conformally transformed function) $\tilde w (x,y)=w(\e^x,y)$.
Moreover, $\tilde w<0$ is defined on a domain that is conformally
close to a half-disk (see lemma 3.2 of \cite{CM8}).  Suppose for a
moment that $\tilde w$ was actually harmonic and defined on the
half-space $\{x>-\log 2\}$, then by the mean value equality (see
fig. \ref{f:f6}) and the sign on $\tilde w$
\begin{equation}  \label{e:expl}
w(\e^x,0)=\tilde w (x,0)=\frac{1}{2\pi x}
\int_{\partial D_x(x,0)}\tilde w
\leq \frac{1}{2\pi x}\int_{D_1\cap \partial D_x(x,0)}\tilde w\, .
\end{equation}
By the Harnack inequality, it would then follow from \eqr{e:expl} that
$w(\e^x,0)=\tilde w(x,0)\leq c_1\,\tilde w(0,0)/x$;
which is the desired lower bound.  The upper bound follows
similarly or by an inversion formula; see section 3 of \cite{CM8}.

\begin{figure}[htbp]
    \setlength{\captionindent}{20pt}
    \begin{minipage}[t]{0.5\textwidth}
    \centering\input{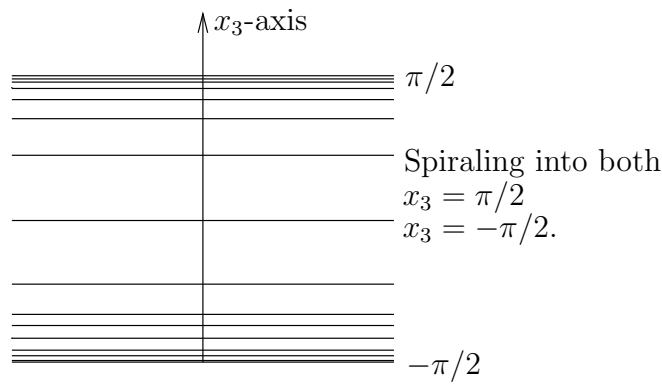}
    \caption{Properness; e.g., need to rule out
    that one of the multi-valued graphs  can
    contain a graph  like  $\arctan (\theta / \log \rho)$,
    where $(\rho , \theta)$ are polar coordinates.} \label{f:f7}
    \end{minipage}
\end{figure}

\vskip2mm
By Theorem \ref{t:main2},
the fastest possible decay for
$w(\rho,0)/w(1,0)$ is $c_1/\log \rho$.  Hence if
$\tilde w(x+iy)=w(\e^x,y)$, then
the fastest possible decay for
$\tilde w(x,0)/\tilde w(0,0)$ is $c_1/x$ and
as mentioned above, then using \eqr{e:5/6} it can be seen that
$\tilde w$ is
almost harmonic.
This decay is achieved for the harmonic function
$\tilde v (z)=-\text{Re}\, z^{-1}=-x/(x^2+y^2)<0$ (see fig. \ref{f:f7}) and if
\begin{equation}  \label{e:exofu}
u(\rho,\theta)=\int_0^{\theta}v(\rho,y)\,dy
=\int_0^{\theta}\tilde v(\log \rho,y)\,dy
= \int_0^{\theta}\frac{-\log \rho\,dy}{(\log\rho)^2+y^2}
=-\arctan  \frac{\theta}{\log \rho} \, ,
\end{equation}
then the graph of $u$ is an embedded $\infty$-valued
harmonic graph lying in a slab, i.e., $|u|\leq \pi/2$,
and hence in particular is not proper.  Note also that if $u$ is given
by \eqr{e:exofu}, then $u_{\theta}=v$
and $w/v$ is uniformly bounded above and below.
We next want to rule out not only this as an example of one half of an embedded
minimal disk, but more generally any $\infty$-valued minimal graph in
a half-space.

\section{Towards removability of singularities}
\label{s:proper}

As mentioned above, we next want to rule out that any $\infty$-valued graph
which lies in a half-space can be one half of an embedded minimal disk.
This is done in Theorem \ref{c:main} below
where the short
curves $\sigma_{\theta}$ in the minimal disk
will connect two multi-valued sub-graphs and thus
 each sub-graph is essentially one half of the disk.
These curves are needed to conclude properness since, as mentioned
above, just having one $\infty$-valued minimal graph in a slab is
not in itself a contradiction.

There is a second key assumption needed to prove this properness:
The outer radii $R_i$ must be going to infinity.  We refer to this
as the global case.  The alternative, the local case, is when the
$R_i$ are bounded. The difference between these two cases is that in
the local case it is possible to have non-proper
limits which spiral infinitely into a plane.  We
will return to this in Part \ref{p:local}.

\begin{figure}[htbp]
    \setlength{\captionindent}{20pt}
    \begin{minipage}[t]{0.5\textwidth}
    \centering\input{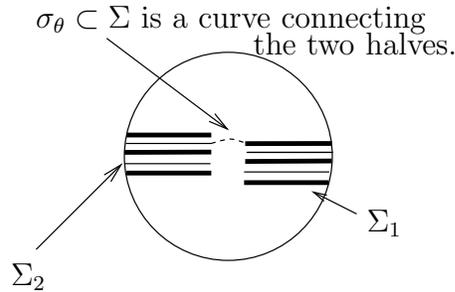}
    \caption{The short curves $\sigma_{\theta}$ in
Theorem \ref{c:main} connecting the two multi-valued graphs.}\label{f:f8}
    \end{minipage}
\end{figure}

Let us illustrate in an example how these curves $\sigma$
could be chosen.  Here, to be consistent with \cite{CM8}, we
use the helicoid which spirals downward so $w<0$.
If $\Sigma$ is the helicoid, i.e.,
$\Sigma=(s\,\cos t,s\,\sin t, -t)$
where $s,\,t\in\RR$, then $\Sigma\setminus \{s=0\}$ consists of
two $\infty$-valued
graphs $\Sigma_1$, $\Sigma_2$ and the curves $\sigma_t=\Sigma\cap \{x_3=t\}$
are short curves connecting the two halves; see fig. \ref{f:f8}.

\begin{Thm}  \label{c:main}
(See corollary 1.2 in \cite{CM8} for the precise statement).
Let $\Sigma_i$ be a sequence as in Theorem \ref{t:t0.1} and suppose
that $\Sigma_{i,1}$, $\Sigma_{i,2}$ are multi-valued graphs
in $\Sigma_i$ that spiral together (one inside the other).  We
claim that if $\Sigma_{i,1}$, $\Sigma_{i,2}$ can be connected by short curves
in $\Sigma_i$ (see fig. \ref{f:f8}), then they cannot spiral into a plane.
That is, they cannot accumulate in finite height.
\end{Thm}

\vskip2mm \noindent {\bf{The idea of the proof of Theorem
\ref{c:main}}}: Let $\Sigma_1$, $\Sigma_2$ be two $\infty$-valued
graphs of $u_1$, $u_2$ that spiral together, are part of an
embedded minimal disk, and can be connected by short curves in the
disk.  We claim that the two graphs must grow out of any
half-space.   Suppose that they are contained in the half-space
$u_i\geq 0$ and that they spiral downward, i.e., $w_i<0$; we will
get a contradiction.

\begin{figure}[htbp]
    \setlength{\captionindent}{20pt}
\begin{minipage}[t]{0.5\textwidth}
    \centering\input{dis9.pstex_t}
    \caption{The flux argument in the proof of Theorem \ref{c:main}:  The
initial height and separation.}
\label{f:flux1}
    \end{minipage}\begin{minipage}[t]{0.5\textwidth}
    \centering\input{dis10.pstex_t}
    \caption{The flux argument in the proof of Theorem \ref{c:main}:  The
curve where the average is calculated over.}
\label{f:flux2}
    \end{minipage}
\end{figure}

Using that the curvature of a multi-valued embedded minimal graph
decays faster than quadratically will allow us (as in Theorem
\ref{t:main2}) to model the graphs $\Sigma_1$, $\Sigma_2$ by
graphs of harmonic functions.  So suppose for a moment that both
$u_1$ and $u_2$ are harmonic.  We will show, using a flux
argument, that if the separation is large compared with the
initial height (see fig. \ref{f:flux1}), then the fact that the
two graphs are part of an embedded minimal disk will eventually
force the graphs to grow out of the half-space that they are
assumed to lie in.  Namely, set
\begin{equation}
    \tilde u_i(x,y)=u_i(\e^x,y)
\end{equation}
 (i.e., make the conformal
change $\rho=\e^x$ and $\theta=y$) and let $Q$ be the
quarter-space $\{ (x,y) \, | \, x , y \geq 0 \}$ and $\tilde D_r =
\{ (x,y)  \, | \, x^2+y^2 \leq r^2 \}$.  To see that each of the
two graphs would grow out of the half-space, set (see fig.
\ref{f:flux2})
\begin{equation}
\text{Av}_i(r)=\frac{1}{r}\int_{\partial \tilde D_r\cap Q}\tilde u_i\, .
\end{equation}
The claim follows once we show that there are constants $C_1$ and $C_2$
depending on the ratio of the initial height with the initial separation
so that for $r\geq 2$
\begin{equation}  \label{e:goalie}
\text{Av}_1(r)+\text{Av}_2(r)
-\text{Av}_1(1)-\text{Av}_2(1)\leq -C_1\,\log^2 r+C_2\,\log r\, .
\end{equation}
Namely, note that if \eqr{e:goalie} holds, then the left hand side
of \eqr{e:goalie} goes to negative infinity as $r \to \infty$,
contradicting that each $\tilde u_i > 0$.  Notice that we can take
$r \to \infty$ to get the contradiction precisely because the
radii $R_i \to \infty$ in Theorem \ref{t:t0.1}. This is the key
place where we use that  $R_i \to \infty$ (see Part \ref{p:local}
for an example which shows this is necessary).

 To see \eqr{e:goalie},
note that when $u_i$ is harmonic, then so is $\tilde u_i$ and by
Stokes' theorem
\begin{equation}  \label{e:I_i'}
\text{Av}_i'(r)=
\frac{1}{r}\int_{\partial \tilde D_r\cap Q}\frac{d\tilde u_i}{dn}
=\frac{1}{r}\int_{0}^{r}\frac{\partial \tilde u_i}{\partial y}(s,0)\,ds
 +\frac{1}{r}\int_{0}^{r}\frac{\partial \tilde u_i}{\partial x}(0,s)\,ds\, .
\end{equation}
Using that the separation is roughly equal to
$2\pi\,(u_i)_{\theta}$ we get by combining \eqr{e:I_i'} with the
lower bound for the separation from Theorem \ref{t:main2} that (more or less)
\begin{equation}    \label{e:h1}
\text{Av}_i'(r)\leq -\frac{C_1\log r}{2\,r}
+\frac{1}{r}\int_{0}^{r}\frac{\partial \tilde u_i}{\partial x}(0,s)\,ds\, .
\end{equation}
Since the two spiraling curves (i.e., $\theta\to (\theta,u_i
(1,\theta))$ for $i=1,2$, see fig. \ref{f:flux3}) together with
the two short curves that are assumed to exist bounds a disk (see
fig. \ref{f:flux4}), we get by Stokes' theorem (see lemma 5.1 of
\cite{CM8}) that
\begin{equation}     \label{e:h2}
\left|  \int_{0}^{r}\frac{\partial \tilde u_1}{\partial
x}(0,s)\,ds +   \int_{0}^{r}\frac{\partial \tilde u_2}{\partial
x}(0,s)\,ds \right| = \left|  \int_{0}^{r}\frac{\partial
u_1}{\partial \rho}(1,\theta)\,d\theta +
\int_{0}^{r}\frac{\partial u_2}{\partial \rho}(1,\theta)\,d\theta
\right| \leq C_2 \, .
\end{equation}
Here, we bounded the flux along the short curves by the length of
the curves since each $|\nabla u_i| \leq 1$ (as is the case for
the restrictions of the coordinate functions).
 Adding the bounds on $\text{Av}_1'$ and $\text{Av}_2'$ in \eqr{e:h1} and
substituting \eqr{e:h2} gives
\begin{equation}  \label{e:sumI}
\text{Av}_1'(r)+\text{Av}_2'(r) \leq -\frac{C_1\log
r}{r}+\frac{C_2}{r} \, .
\end{equation}
Integrating \eqr{e:sumI} gives \eqr{e:goalie}.

\begin{figure}[htbp]
    \setlength{\captionindent}{20pt}
\begin{minipage}[t]{0.5\textwidth}
    \centering\input{dis11.pstex_t}
    \caption{The flux argument in the proof of Theorem \ref{c:main}:
One of the two spiraling curves; $\theta\to (\theta,u_i(1,\theta))$.}
\label{f:flux3}
    \end{minipage}\begin{minipage}[t]{0.5\textwidth}
    \centering\input{dis12.pstex_t}
    \caption{The flux argument in the proof of Theorem \ref{c:main}:  The
two short curves together with the two spiraling curves bounds a
disk in $\Sigma$.}
\label{f:flux4}
    \end{minipage}
\end{figure}

In the general case where $u_i$ satisfies the minimal surface equation there
are a number of difficulties that have to be dealt with;
see lemma 4.1 of \cite{CM8}.

\vskip2mm
Now this was a little analysis of multi-valued solutions of the
minimal surface equation and can all be found in \cite{CM8}.  It
illustrates that once we show existence of multi-valued minimal
graphs in embedded minimal disks and know that such graphs extend as
graphs with a sufficiently rapidly growing number of sheets, then we
get a removable
singularity theorem.

\section{Existence of multi-valued graphs and the one-sided curvature
estimate} \label{s:13}

We now come to our key results for embedded minimal disks. These
are some of the main ingredients in the proof of Theorem
\ref{t:t0.1}.
 The
first says that if the curvature of such a disk $\Sigma$ is large
at some point $x\in \Sigma$, then nearby $x$
 a multi-valued graph forms (in $\Sigma$) and this extends
(in $\Sigma$) almost all the way to the boundary.  Precisely this is:

\begin{Thm} \label{t:blowupwinding0}
(Theorem $0.2$ in \cite{CM4}).  See fig. \ref{f:f9} and fig. \ref{f:f10}.
Given $N\in \ZZ_+$, $\epsilon > 0$, there exist
$C_1,\,C_2>0$ so: Let
$0\in \Sigma \subset B_{R}\subset \RR^3$ be an embedded minimal
disk, $\partial \Sigma\subset \partial B_{R}$. If
$\max_{B_{r_0} \cap \Sigma}|A|^2\geq 4\,C_1^2\,r_0^{-2}$ for some
$R>r_0>0$, then there exists
(after a rotation)
an $N$-valued graph $\Sigma_g \subset \Sigma$ over $D_{R/C_2}
\setminus D_{2r_0}$ with gradient $\leq \epsilon$
and
$\Sigma_g \subset \{ x_3^2 \leq \epsilon^2 \, (x_1^2 + x_2^2) \}$.
\end{Thm}

\begin{figure}[htbp]
    \setlength{\captionindent}{20pt}
    \begin{minipage}[t]{0.5\textwidth}
\centering\input{shn3.pstex_t}
    \caption{Part 1 of the proof of Theorem \ref{t:blowupwinding0};
see Theorem \ref{t:blowupwindinga} - finding a small
    multi-valued graph in a disk near a point of large curvature.}\label{f:f9}
\end{minipage}\begin{minipage}[t]{0.5\textwidth}
    \centering\input{shn2.pstex_t}
    \caption{Part 2 of the proof of Theorem \ref{t:blowupwinding0}; see
Theorem \ref{t:spin4ever2} - extending a small multi-valued graph
    in a disk.}\label{f:f10}
    \end{minipage}
\end{figure}

As a consequence of Theorem \ref{t:blowupwinding0}, one easily
gets that if $|A|^2$ is blowing up near $0$ for a sequence of
embedded minimal disks $\Sigma_i$, then there is a sequence of
$2$-valued graphs $\Sigma_{i,d}\subset \Sigma_i$, where the
$2$-valued graphs start off on a smaller and smaller scale
(namely, $r_0$ in Theorem \ref{t:blowupwinding0} can be taken to
be smaller as the curvature gets larger). Consequently, by the
sublinear separation growth, such $2$-valued graphs collapse and,
hence, a subsequence converges to a smooth minimal graph through
$0$.  To be precise, given any fixed $\rho > 0$, \eqr{e:slg}
bounds the separation $w$ at $(\rho , 0)$ by
\begin{equation}    \label{e:closeup}
    |w(\rho , 0)| \leq \left( \frac{\rho}{r_0} \right)^{\alpha} \,
    |w(r_0,0)| \leq 2 \pi \, \epsilon \, \rho^{\alpha} \,
    r_0^{1-\alpha} \, ,
\end{equation}
and this goes to $0$ as $r_0 \to 0$ since $\alpha < 1$.  The bound
$|w(r_0,0)| \leq 2 \pi \, \epsilon \, r_0$ in \eqr{e:closeup} came
from integrating the gradient bound on the graph around the circle
of radius $r_0$.  (Here $0$ is a removable singularity for the
limit.) Moreover, if the sequence of such disks is as in Theorem
\ref{t:t0.1}, i.e., if $R_i\to \infty$, then the minimal graph in
the limit is entire and hence, by Bernstein's theorem (theorem
$1.16$ in \cite{CM1}), is a plane.

\begin{figure}[htbp]
    \setlength{\captionindent}{4pt}
    \begin{minipage}[t]{0.5\textwidth}
    \centering\input{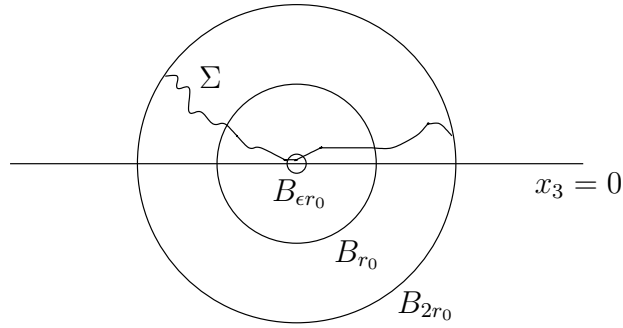}  
    \caption{Theorem \ref{t:t2} - the one-sided curvature estimate for an
embedded minimal disk $\Sigma$ in a half-space with
$\partial \Sigma\subset \partial B_{2r_0}$:  The components of
$B_{r_0}\cap \Sigma$ intersecting $B_{\epsilon r_0}$ are graphs.}\label{f:f11}
    \end{minipage}
\end{figure}

\vskip2mm The second key result is the  curvature estimate for
embedded minimal disks in a half-space.  This theorem says roughly that if
an embedded minimal disk lies in a half-space above a plane and
comes close to the plane, then it is a graph over the plane.
Precisely, this is the following theorem:

\begin{Thm}  \label{t:t2}
(Theorem 0.2 in \cite{CM6}). See fig. \ref{f:f11}.
There exists $\epsilon>0$, such that if
$\Sigma \subset B_{2r_0} \cap \{x_3>0\}
\subset \RR^3$ is an
embedded minimal
disk with $\partial \Sigma\subset \partial B_{2 r_0}$,
then for all components $\Sigma'$ of
$B_{r_0} \cap \Sigma$ which intersect $B_{\epsilon r_0}$
\begin{equation}  \label{e:graph}
\sup_{\Sigma'} |A_{\Sigma}|^2
\leq r_0^{-2} \, .
\end{equation}
\end{Thm}

Using the minimal surface equation and that $\Sigma'$ has points
close to a plane, it is not hard to see that, for $\epsilon>0$
sufficiently small, \eqr{e:graph} is equivalent to the statement
that $\Sigma'$
is a graph over the plane $\{x_3=0\}$.

An embedded minimal surface $\Sigma$ which is as in
Theorem \ref{t:t2} is said to satisfy the
$(\epsilon , r_0)$-{\it effective one-sided Reifenberg condition}; cf.
appendix A of \cite{CM6} and the appendix of \cite{ChC}.
We will often refer to Theorem \ref{t:t2}
as {\it the one-sided curvature estimate}.

\begin{figure}[htbp]
    \setlength{\captionindent}{20pt}
    \begin{minipage}[t]{0.5\textwidth}
    \centering\input{pl2a.pstex_t}
    \caption{The catenoid given by revolving $x_1= \cosh x_3$
around the $x_3$-axis.}  \label{f:f12}
    \end{minipage}\begin{minipage}[t]{0.5\textwidth}
    \centering\input{unot7.pstex_t}
    \caption{Rescaling the catenoid shows that simply connected
(and embedded) is
    needed in the one-sided curvature estimate.}  \label{f:f13}
    \end{minipage}%

\end{figure}

Note that the assumption in Theorem \ref{t:t2}
that $\Sigma$ is simply connected is crucial
as can be seen from the example of a rescaled catenoid. The catenoid
is the minimal surface in $\RR^3$ given by
$(\cosh s\, \cos t,\cosh s\, \sin t,s)$
where $s,t\in\RR$; see fig. \ref{f:f12}.
Under rescalings this converges (with multiplicity two) to
the flat plane; see fig. \ref{f:f13}. Likewise, by considering
the universal cover of the catenoid, one sees that embedded,
and not just immersed, is needed in Theorem \ref{t:t2}.

As an almost immediate consequence of Theorem \ref{t:t2} and a simple
barrier argument we get that
if in a ball two embedded minimal disks come close to each other near
the center of the ball then each of the disks are graphs.  Precisely,
this is the following:

\begin{figure}[htbp]
    \setlength{\captionindent}{20pt}
    \begin{minipage}[t]{0.5\textwidth}
    \centering\input{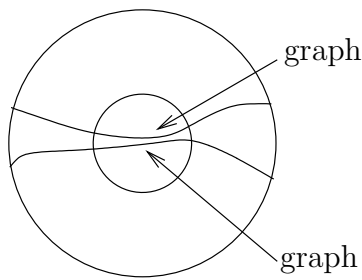}
    \caption{Corollary \ref{c:barrier}:  Two sufficiently close
    components of an embedded minimal disk must each be a graph.}\label{f:f14}
    \end{minipage}
\end{figure}

\begin{Cor}         \label{c:barrier}
(Corollary 0.4 in \cite{CM6}). See fig. \ref{f:f14}.
There exist $c > 1$, $\epsilon >0$ so:
Let $\Sigma_1,\, \Sigma_2 \subset B_{cr_0} \subset \RR^3$ be
disjoint embedded minimal
surfaces with $\partial \Sigma_i \subset \partial B_{cr_0}$
and $B_{\epsilon \, r_0 } \cap \Sigma_i \ne \emptyset$.
If $\Sigma_1 $ is a disk,
 then for all components $\Sigma_1'$ of
$B_{r_0} \cap \Sigma_1$ which intersect
$B_{\epsilon \, r_0}$
\begin{equation}        \label{e:onece}
    \sup_{\Sigma_1'}   |A|^2
        \leq  r_0^{-2}  \, .
\end{equation}
\end{Cor}

Theorem \ref{t:t2} is used to show that the points of large curvature
in an embedded minimal disk all lie on a Lipschitz curve.  To be able
to discuss this and explain why this follows from the theorem lets
introduce some notation for cones.

If $\delta>0$ and $z\in \RR^3$,
then we denote by
$\cone_{\delta}(z)$ the (convex) cone with vertex $z$, cone angle
$(\pi/2 - \arctan \delta)$, and axis parallel to the $x_3$-axis.
That is, see fig. \ref{f:f15a},
\begin{equation}
\cone_{\delta}(z)=\{x\in \RR^3\,|\,x_3^2 \geq
\delta^2\,(x_1^2+x_2^2) \} +z\, .
\end{equation}

\vskip2mm
The next direct consequence of Theorem \ref{t:t2}
(with $\Sigma_d$ playing the role of the plane $x_3 = 0$) will be
needed when we sketch the proof of Theorem \ref{t:t0.1} in the
next section.  This consequence says that the points of large
curvature in an embedded minimal disk have what we will in the next
section call the cone property.  Namely, given a point of large
curvature in such a disk, the next corollary asset that all the other
points of large curvature lies within a double convex cone with vertex
at the initial point of large curvature.  This is the result that will
eventually give the regularity of the singular set (the set of points
of large curvature).  Precisely, this consequence of  Theorem
\ref{t:t2} is the following:

\begin{figure}[htbp]
    \setlength{\captionindent}{20pt}
    \begin{minipage}[t]{0.5\textwidth}
    \centering\input{dis19.pstex_t}
    \caption{The cone $\cone_{\delta}(x)$.}  \label{f:f15a}
    \end{minipage}\begin{minipage}[t]{0.5\textwidth}
    \centering\input{dis20.pstex_t}
    \caption{Corollary \ref{c:conecor}:  With $\Sigma_d$ playing
    the role of $x_3=0$,
    by the one-sided estimate, $\Sigma$ consists of multi-valued
    graphs away from a cone.}   \label{f:f16}
    \end{minipage}
\end{figure}

\begin{Cor}   \label{c:conecor}
(Corollary I.1.9 in \cite{CM6}).  See fig. \ref{f:f16}.
There exists $\delta_0>0$ so: Suppose $\Sigma\subset B_{2R}$,
$\partial \Sigma\subset \partial B_{2R}$ is an embedded minimal
disk containing a $2$-valued graph $\Sigma_d \subset \{x_3^2 \leq
\delta_0^2\, (x_1^2+x_2^2)\}$ over $D_{R}\setminus D_{r_0}$ with
gradient $\leq \delta_0$.  Then each component of $B_{R/2}\cap
\Sigma\setminus (\cone_{\delta_0}(0)\cup B_{2 r_0})$ is a
multi-valued graph with gradient $\leq 1$.
\end{Cor}

Fig. \ref{f:f16} illustrates how this corollary follows from
Theorem \ref{t:t2}.  In this picture, $B_s(y)$ is a ball away from
$0$ and $\Sigma'$ is a component of $B_s(y)\cap \Sigma$ disjoint
from $\Sigma_d$.  It follows easily from the maximum principle
that $\Sigma'$ is topologically a disk.  Since $\Sigma'$ is
assumed to contain points near $\Sigma_d$, then we can let a
component of $B_s(y)\cap \Sigma_d$ play the role of the plane
$\{x_3=0\}$ in Theorem \ref{t:t2} and the corollary follows.

Note that, since $\Sigma$ is compact and embedded, the multi-valued
graphs given by Corollary \ref{c:conecor} spiral through the cone.
Namely, if a graph did close up, then the graph containing
$\Sigma_d$ would be forced to accumulate into it, contradicting
compactness.

\vskip2mm
It can also be seen from Corollary \ref{c:conecor} (see corollary 6.3 of
\cite{CM8}) that
if $\Sigma_d$ and $\Sigma$ are as in Corollary \ref{c:conecor}, then
$\Sigma_d$ extends in $\Sigma \setminus \cone_{\delta_0}(0)$
to a multi-valued graph with at
least $(\log\rho)^2$ many sheets.  Thus Theorem \ref{t:main2} applies.

\section{Regularity of the singular set and Theorem \ref{t:t0.1}}
\label{s:ref1}

In this section we will indicate how to define the singular set $\cS$
in Theorem \ref{t:t0.1} and show the regularity of $\cS$.

First, by a very general compactness argument, we have that for any
sequence of surfaces in $\RR^3$ (minimal or not),
after possibly going to a subsequence,
then there is a well defined notion of points where the second
fundamental form of the sequence blows up.

\begin{Lem} \label{l:inftyornot}
Let $\Sigma_i\subset B_{R_i}$, $\partial \Sigma_i\subset \partial B_{R_i}$,
and $R_i\to \infty$ be a sequence of (smooth) compact surfaces.
After passing to a subsequence, $\Sigma_j$, we
may assume that for each $x\in \RR^3$ either (a) or (b) holds:\\
(a) $\sup_{B_{r}(x)\cap
\Sigma_j}|A|^2\to \infty$ for all $r>0$,\\
(b) $\sup_j\sup_{B_r(x)\cap
\Sigma_j}|A|^2<\infty$ for some $r>0$.
\end{Lem}

\begin{proof}
For $r>0$ and an integer $n$, define a sequence of functions on $\RR^3$ by
\begin{equation}
\cA_{i,r,n}(x)=\min \{ n,\sup_{B_r(x)\cap \Sigma_i}|A|^2\}\, ,
\end{equation}
where we set $\sup_{B_r(x)\cap \Sigma_i}|A|^2=0$ if
$B_r(x)\cap \Sigma_i=\emptyset$.
Set
\begin{equation}
\cD_{i,r,n}=\lim_{k\to \infty} 2^{-k} \,
\sum_{m=0}^{2^k-1} \cA_{i,(1+m2^{-k})r,n}\, ,
\end{equation}
then $\cD_{i,r,n}$ is continuous and $\cA_{i,2r,n}\geq \cD_{i,r,n}\geq
\cA_{i,r,n}$.
 Let $\nu_{i,r,n}$
be the (bounded) functionals,
\begin{equation}
\nu_{i,r,n}(\phi)=\int_{B_n}\phi\,\cD_{i,r,n}\text{ for }
\phi\in L^2 (\RR^3)\, .
\end{equation}
By standard compactness for fixed $r$, $n$, after passing to a subsequence,
$\nu_{j,r,n}\to \nu_{r,n}\text{ weakly}$.
In fact (since the unit ball in $L^2(\RR^3)$ has a countable basis),
by an easy diagonal argument, after passing to a
subsequence, we may assume that for all $n, m \geq 1$ fixed
$\nu_{j,2^{-m},n}\to \nu_{2^{-m},n}\text{ weakly}$.
Note that if $x\in \RR^3$ and for all $m$, $n$ with $n\geq |x|+1$,
(identify $B_{2^{-m}}(x)$ with its characteristic function)
\begin{equation} \label{e:inftyornot}
\nu_{2^{-m},n}(B_{2^{-m}}(x))\geq n\, \Vol (B_{2^{-m}})\, ,
\end{equation}
then for each fixed $r>0$, $\sup_{B_r(x)\cap\Sigma_j}|A|^2\to \infty$.
Conversely, if for some $n \geq |x| +1$, $m$, \eqr{e:inftyornot}
fails at $x$, then
$\sup_j\sup_{B_r(x)\cap
\Sigma_j}|A|^2<\infty$ for $r=2^{-m-1}$.
\end{proof}

From this lemma we have that after possibly passing to a subsequence,
then there is a well defined notion of the set of points where the
curvatures blow up of a given sequence of embedded minimal disks.  To
show that this set is in fact a Lipschitz curve we will see below
that,
as a consequence of the one sided curvature estimate, the set of
such points has what we will call the cone
property.

Fix $\delta>0$. We will say that a subset $\cS\subset \RR^3$ has the
cone property (or the $\delta$-cone property) if $\cS$ is closed and
nonempty and the following holds:\\
(1) If $z\in \cS$, then $\cS\subset
\cone_{\delta}(z)$.\\
(2) If $t\in x_3(\cS)$ and
$\epsilon>0$, then $\cS\cap \{t<x_3<t+\epsilon\}\ne \emptyset$ and
$\cS\cap \{t-\epsilon<x_3<t\}\ne \emptyset$.

\vskip2mm
Note that (2) just says that each point in $\cS$ is the limit of points
coming from above and below.

When $\Sigma_i\subset B_{R_i}\subset \RR^3$  is a sequence of embedded
minimal disks with $\partial \Sigma\subset \partial B_{R_i}$,
$R_i\to \infty$ and $\Sigma_j$ is the subsequence given by Lemma
\ref{l:inftyornot}
and $\cS$ is the sets of points where the curvatures of $\Sigma_j$
blow up (i.e., where (a) in Lemma
\ref{l:inftyornot} holds), then (as we indicated above)
we will see below that  $\cS$ has the
cone property (after a rotation of $\RR^3$).   Hence (by the next lemma),
$\cS$ is in this case a Lipschitz curve which
is a graph over the $x_3$-axis.  Note that in the case where $\Sigma_i$
is a sequence of rescaled helicoids, then $\cS$ is simply the $x_3$-axis.

\begin{figure}[htbp]
    \setlength{\captionindent}{20pt}
    \begin{minipage}[t]{0.5\textwidth}
    \centering\input{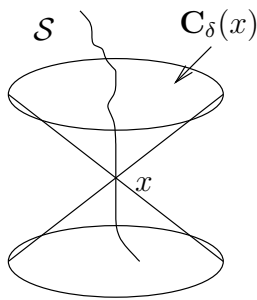}
    \caption{It follows from the one-sided curvature estimate that the singular
    set has the cone property and hence is a Lipschitz curve.}  \label{f:f15}
    \end{minipage}
\end{figure}

\begin{Lem}  \label{l:regsing}
See fig. \ref{f:f15}.
If $\cS\subset \RR^3$ has the $\delta$-cone property, then
$\cS\cap
\{x_3=t\}$ consists of exactly one point $\cS_t$ for all
$t\in\RR$, and $t\to \cS_t$ is a Lipschitz parameterization of
$\cS$ with
\begin{equation}  \label{e:lipscbd}
|t_2-t_1|\leq |\cS_{t_2}-\cS_{t_1}|\leq \sqrt{1+\delta^{-2}}\,|t_2-t_1|\, .
\end{equation}
\end{Lem}

\begin{proof}
Since $\cS$ is nonempty, we may after translation assume that $0\in\cS$.
By (1), it follows that $\cS\cap \{x_3=t\}$
consists of at most one point for each $t\in \RR$.  Assume that
$\cS\cap \{x_3=t_0\}=\emptyset$ for some $t_0$.  Since $\cS\subset
\RR^3$ is a nonempty closed set and $x_3:\cS\subset
\cone_{\delta}(0)\to \RR$ is proper, then $x_3(\cS)\subset \RR$ is
also closed (and nonempty). Let $t_s\in x_3(\cS)$ be the closest
point in $x_3(\cS)$ to $t_0$. The desired contradiction now easily
follows since either $\cS\cap \{t_s<x_3<t_0\}$ or $\cS\cap
\{t_0<x_3<t_s\}$ is nonempty by assumption.

It follows that $t\to \cS_t$ is a well-defined curve (from $\RR$
to $\cS$). Moreover, since
\begin{equation}
    \cS_{t_2}\subset
    \{x_3=t_1+(t_2-t_1)\}\cap \cone_{\delta}(\cS_{t_1}) \subset
    B_{\sqrt{1+\delta^{-2}}|t_2-t_1|}(\cS_{t_1}) \, ,
\end{equation}
\eqr{e:lipscbd} follows.
\end{proof}

Suppose next that $\Sigma_i$ is as in Theorem \ref{t:t0.1}, that
is, $\Sigma_i \subset B_{R_i}=B_{R_i}(0)\subset \RR^3$ is a
sequence of embedded minimal disks with $\partial \Sigma_i\subset
\partial B_{R_i}$ where $R_i\to \infty$ and $\sup_{B_1\cap
\Sigma_i}|A|^2\to \infty$.  Let $\Sigma_j$ and $\cS$ be the
subsequence and set, respectively,
 given by Lemma \ref{l:inftyornot}
($\cS$ is the set of points where (a) holds in Lemma
\ref{l:inftyornot}).  In particular, $\cS$ is closed by definition
and nonempty by the assumption of Theorem \ref{t:t0.1}. From
Corollary \ref{c:conecor} (cf. the sketch of the proof of Theorem
\ref{t:t0.1} below), it follows that (1) above holds, so to see
that $\cS$ has the cone property all we need to see is that (2)
holds.  This follows from the next lemma which relies in part on
Theorem \ref{c:main}:

\begin{figure}[htbp]
    \setlength{\captionindent}{20pt}
    \begin{minipage}[t]{0.5\textwidth}
\centering\input{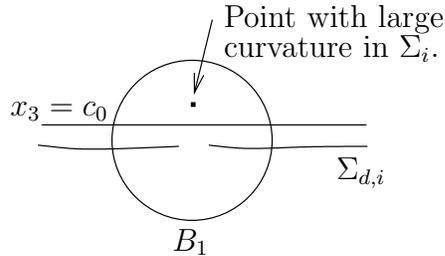}
    \caption{Lemma \ref{l:blowupnear} - points with large
curvature in $\Sigma_i$
above the plane $x_3=c_0$ but near the center of the $2$-valued
graphs $\Sigma_{d,i}$.}
\label{f:f17}
\end{minipage}
\end{figure}

\begin{Lem}  \label{l:blowupnear}
(Lemma I.1.10 in \cite{CM6}).  See fig. \ref{f:f17}.
There exists $c_0>0$ so:
Let $\Sigma_i\subset B_{R_i}$, $\partial \Sigma_i\subset \partial B_{R_i}$
be a sequence of embedded minimal disks with $R_i\to \infty$.
If $\Sigma_{d,i}\subset \Sigma_i$ is a sequence of $2$-valued
graphs over $D_{R_i/C}\setminus D_{\epsilon_i}$ with $\epsilon_i\to 0$
and $\Sigma_{d,i}\to \{x_3=0\}\setminus \{0\}$, then
\begin{equation}
\sup_{B_1 \cap \Sigma_i\cap \{x_3>c_0 \}}|A|^2\to \infty\, .
\end{equation}
\end{Lem}

\vskip2mm \noindent {\bf{Sketch of the proof of Theorem
\ref{t:t0.1}}}: From all of these results above, we know that if
$\Sigma_i$ is a sequence as in Theorem \ref{t:t0.1} and
$\Sigma_j$, $\cS$ are given by Lemma \ref{l:inftyornot}, then
$\cS$ has the cone property and hence, by Lemma
\ref{l:inftyornot}, $\cS$ is a Lipschitz graph over the
$x_3$-axis.  As mentioned above, it also follows from Theorem
\ref{t:blowupwinding0} together with Bernstein's theorem that, for
each $x\in \cS$ and each $j$ sufficiently large, there is a
$2$-valued graph in $\Sigma_j$ and that this sequence of
$2$-valued graphs converges (after possibly going to a
subsequence) to a plane through $x$.  Since $\Sigma_j$ is
embedded, all of these planes coming from different points $x\in
\cS$ must be parallel and it now follows from the one sided
curvature estimate that $\Sigma_j$ consists of multi-valued graphs
away from $\cS$ .  It is easy to see that there must be at least
two such multi-valued graphs in the complement of $\cS$.  That
there are not more than two follows from a barrier argument; see
proposition II.1.3 of \cite{CM6}. This completes the rough sketch
of the proof of Theorem \ref{t:t0.1}.

\part{The proof of the existence of multi-valued graphs}

Before we proceed let us briefly review the strategy of the proof of
our main theorem that every embedded
minimal disk is either a graph of
a function or part of a double spiral staircase.

\begin{figure}[htbp]
    \setlength{\captionindent}{20pt}
    \begin{minipage}[t]{0.5\textwidth}
    \centering\input{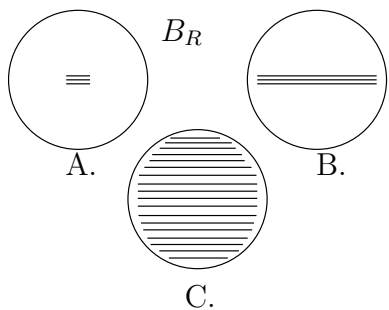}
    \caption{Proving Theorem \ref{t:t0.1}.  A. Finding a small
$N$-valued graph in $\Sigma$.
    B. Extending it in $\Sigma$ to a large
    $N$-valued graph.  C. Extending
    the number of sheets.}\label{f:f3a}
    \end{minipage}
\end{figure}

The proof has the following three main steps; see fig.
\ref{f:f3a}:
\vskip1.5mm
\noindent
A.  Fix an integer $N$ (the
``large'' of the curvature in what follows will depend on $N$).
If an embedded minimal disk $\Sigma$ is not a graph (or
equivalently if the curvature is large at some point), then it
contains an $N$-valued minimal graph which initially is shown to
exist on the scale of $1/\max |A|$.  That is, the $N$-valued graph
is initially shown to be defined on an annulus with both inner and
outer radius inversely proportional to $\max |A|$.
\vskip1.5mm
\noindent
B.  Such a potentially small $N$-valued graph sitting
inside $\Sigma$ can then be seen to extend as an $N$-valued graph
inside $\Sigma$ almost all the way to the boundary.  That is, the
small $N$-valued graph can be extended to an $N$-valued graph
defined on an annulus where the outer radius of the annulus is
proportional to $R$.  Here $R$ is the radius of the ball in
$\RR^3$ that the boundary of $\Sigma$ is contained in.
\vskip1.5mm
\noindent
C.  The $N$-valued graph not only extends horizontally
(i.e., tangent to the initial sheets) but also vertically (i.e.,
transversally to the sheets).  That is, once there are $N$ sheets
there are many more and, in fact, the disk $\Sigma$ consists of
two multi-valued graphs glued together along an axis. \vskip1.5mm

To describe the existence, i.e., A., of multi-valued graphs in embedded
minimal disks, we will need the notion of a blow up point. Roughly
speaking, a blow up point is a point where the curvature is, up to
a fixed constant, the maximum.

\begin{figure}[htbp]
    \setlength{\captionindent}{20pt}
    \begin{minipage}[t]{0.5\textwidth}
\centering\input{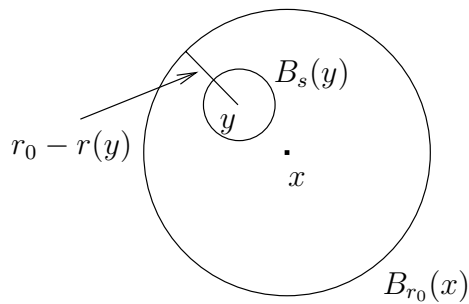}
    \caption{Existence of blow up points, that is pairs of points
$y\in \Sigma$ and $s>0$ satisfying \eqr{e:ti}.}
\label{f:ti}
\end{minipage}
\end{figure}

To make this precise, let $x\in \Sigma\subset B_{r_0}(x)\subset
\RR^3$ be a smooth (compact) surface embedded, or just immersed,
with $\partial \Sigma\subset \partial B_{r_0}(x)$. Here
$B_{r_0}(x)$ is the extrinsic ball of radius $r_0$ but could as
well have been an intrinsic ball in which case the notion of a
blow up point below would have to be appropriately changed.
Suppose that $|A|^2(x)\geq 4 \, C^2\,r_0^{-2}$ for some constant
$C>0$.  We claim that there is a ``blow up point'' $y\in
B_{r_0}(x)\cap \Sigma$ and $s>0$ such that $B_s(y)\subset
B_{r_0}(x)$ and
\begin{equation}  \label{e:ti}
\sup_{B_s(y)\cap \Sigma}|A|^2\leq 4\,C^2\,s^{-2}= 4\,|A|^2(y)\, .
\end{equation}
That is, the curvature at $y$ is large (this just means that $C$ should
be thought of as a large constant) and is almost
(up to the constant $4$) the maximum on
the ball $B_s(y)$.

\begin{proof} (of claim) That there exists such a point $y$ is easy to see; on
$B_{r_0}(x)\cap \Sigma$ set $F(z)=(r_0-r(z))^2\,|A|^2(z)$ where
$r(z)=|z-x|$.  Then
\begin{equation}
F(x)\geq 4\,C^2\, , F\geq 0\, , \text{ and }
F\left. \right|_{\partial B_{r_0}(x)\cap \Sigma}=0\, .
\end{equation}
Let $y$ be where the maximum
of $F$ is achieved and set $s= C/|A|(y)$.  One easily checks that
$y$, $s$ have the required properties; see fig. \ref{f:ti}.
Namely, clearly $|A|^2(y)= C^2\,s^{-2}$
and since $y$ is where the maximum of $F$ is achieved,
\begin{equation}
|A|^2(z)\leq
\left(\frac{r_0-r(y)}{r_0-r(z)}\right)^2\,|A|^2(y)\, .
\end{equation}
Since $F(x)\geq 4\,C^2$ it follows from
the choice of $s$ that $|r_0-r(y)|\leq 2\,|r_0-r(z)|$ for
$z\in B_s(y)\cap \Sigma$.  Hence,
$|A|^2(z)\leq 4\,|A|^2(y)$.  Together this gives \eqr{e:ti}.
\end{proof}

Returning to the existence of multi-valued graphs, then the way we
showed Theorem \ref{t:blowupwinding0} was by combining a blow up
result with the extension of multi-valued graphs proven in
\cite{CM3}. This blow up result says that if an embedded minimal
disk in a ball has large curvature at a point, then it contains a
small (in fact on the scale of one over the square root of the
curvature) almost flat multi-valued graph nearby, that is:

\begin{Thm} \label{t:blowupwindinga}
(Theorem $0.4$ in \cite{CM4}).  See fig. \ref{f:f9}.
Given $N , \omega>1$, and $\epsilon > 0$, there exists
$C=C(N,\omega,\epsilon)>0$ so: Let
$0\in \Sigma \subset B_{R}\subset \RR^3$ be an embedded minimal
disk, $\partial \Sigma\subset \partial B_{R}$. If
$\sup_{B_{r_0} \cap \Sigma}|A|^2\leq 4\,C^2\,r_0^{-2}=4\,|A|^2(0)$
for some $0<r_0<R$, then there exist
$ \bar{R} < r_0 / \omega$ and (after a rotation)
an $N$-valued graph $\Sigma_g \subset \Sigma$ over $D_{\omega \bar{R} }
\setminus D_{\bar{R} }$ with gradient $\leq \epsilon$, and
$\dist_{\Sigma}(0,\Sigma_g) \leq 4 \, \bar{R}$.
\end{Thm}

Note that $0$ in Theorem \ref{t:blowupwindinga} is a
blow up point in the above sense.

\vskip2mm
 The result that we needed from \cite{CM3} (combining
theorem 0.3 and lemma II.3.8 there) is Theorem \ref{t:spin4ever2}
below that allows us to extend the (small) multi-valued graphs
given by Theorem \ref{t:blowupwindinga} almost out to the boundary
of the ``big'' ball $B_R$. In this theorem, by the middle sheet
$(\Sigma_g)^M$ of an $N$-valued graph $\Sigma_g$ we mean the
portion over
\begin{equation}
\{ (\rho ,\theta ) \in \cP \, | \, r_1 < \rho < r_2 {\text{ and }}
0 \leq \theta \leq 2 \, \pi \} \, .
\end{equation}

\begin{Thm} \label{t:spin4ever2}
\cite{CM3}.  See fig. \ref{f:f10}.
Given $N_1$ and $\tau > 0$, there exist $N , \Omega, \epsilon > 0$ so:
If $\Omega \, r_0 < 1 < R / \Omega$, $\Sigma \subset B_{R}$ is
an embedded minimal disk with $\partial \Sigma \subset \partial B_{R}$,
and $\Sigma$ contains an $N$-valued minimal graph $\Sigma_g$ over
$D_1 \setminus D_{r_0}$ with gradient $\leq \epsilon$ and
$\Sigma_g \subset \{ x_3^2 \leq \epsilon^2 (x_1^2 + x_2^2) \}$, then
$\Sigma$ contains a $N_1$-valued graph $\Sigma_d$ over
$D_{R/\Omega} \setminus D_{r_0}$ with gradient $\leq \tau$ and
$(\Sigma_g)^M \subset \Sigma_d$.
\end{Thm}

\section{The proof of existence of small multi--valued graphs -
Theorem \ref{t:blowupwindinga}}
\label{s:blowup}

We will here describe some of the ideas that go into the proof of
the existence of the small multi-valued graphs near a point of
large curvature.  That is, the  proof of Theorem
\ref{t:blowupwindinga}.  In this theorem, $\Sigma$ is an embedded
minimal disk and $0 \in \Sigma$ is a blow up point with scale
$r_0$ satisfying
\begin{equation}    \label{e:bbb}
    \sup_{B_{r_0} \cap \Sigma}|A|^2 \leq
        4\,C^2\,r_0^{-2}=4\,|A|^2(0) \, .
\end{equation}
The theorem then gives a multi-valued graph $\Sigma_g$ contained
in $\Sigma$ defined  on the scale $r_0$ (\eqr{e:bbb} says that
$r_0$ is proportional to one over the squareroot of the curvature)

\begin{figure}[htbp]
    \setlength{\captionindent}{20pt}
    \begin{minipage}[t]{0.5\textwidth}
\centering\input{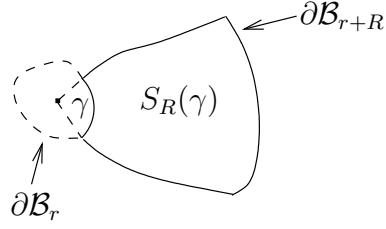}
    \caption{An intrinsic sector over a curve $\gamma$.}
\label{f:9}
\end{minipage}
\end{figure}

The key step in finding the multi-valued graph is to find many
large pieces of $\Sigma$ with a (scale-invariant) quadratic
curvature bound (these pieces will be intrinsic sectors; see fig.
\ref{f:9}). To do this,  we use the upper bound on $|A|^2$ in
\eqr{e:bbb} to prove that the area of intrinsic balls in $\Sigma$
grows polynomially and, consequently, get an average curvature
bound. This average curvature bound and a curvature estimate for
embedded disks (see Proposition \ref{p:scsi} below) will give
large pieces of $\Sigma$ with the desired quadratic curvature
bound. Using the lower bound on $|A|^2 (0)$ in \eqr{e:bbb}, we
show that there are many such pieces so that two must be close
together in $\RR^3$; embeddedness implies that these are disjoint,
hence almost stable, and therefore nearly flat. Piecing together
these large flat pieces then gives the desired $N$-valued graph.

Throughout this section, $\Sigma \subset B_R$ is an embedded minimal
disk with $\partial \Sigma\subset \partial B_R$.

\vskip2mm  Before discussing the main steps in the proof, we will
need to recall three facts about minimal surfaces from \cite{CM4}.
The first is the relationship between area and total curvature of
intrinsic balls for disks with nonpositive curvature. The second
is a curvature estimate for embedded minimal disks with bounded
total curvature.   The third is that nearby, but disjoint, minimal
surfaces with bounded curvature must in fact be almost stable.

\vskip2mm \noindent {\underline{Area and total curvature}}. The
relationship between area and total curvature is particularly
simple for disks.  Essentially, this is because the first
variation of length of a geodesic circle is (up to a constant)
given by the total curvature of the disk using the Gauss-Bonnet
theorem; this can also be seen using the Jacobi equation for
geodesics. Namely, as in corollary $1.7$ of \cite{CM4},
integrating the Jacobi equation (for geodesics) and using
$K_{\Sigma} = -|A|^2/2$ gives
\begin{equation}    \label{e:jacga}
 4\,(\Area \, (\cB_{R}) - \pi \, R^2) = 2
\int_{0}^{R} \int_0^{t} \int_{\cB_{s}} |A|^2
    = \int_{\cB_{R}}|A|^2\,(R-r)^2 \, ,
\end{equation}
where $r(x) = \dist_{\Sigma}(0,x)$.  The second equality in
\eqr{e:jacga} used two integrations by parts (i.e., $\int_0^R f(t)
\, g''(t) \, dt = \int_0^R f''(t) \, g(t) \, dt$ with $f(t) =
\int_0^t \int_{\cB_s} |A|^2$ and $g(t) = (R-t)^2$).

We will see that \eqr{e:jacga} often leads to bounds on the area
of the ball $\cB_R$.
 For example, when $\cB_R$ is stable, then using $R-r$ (which vanishes on
 $\partial \cB_R$)
 in the stability inequality (see \eqr{e:stabin0} below) gives
\begin{equation}    \label{e:jacgab}
 4\,(\Area \, (\cB_{R}) - \pi \, R^2)
    = \int_{\cB_{R}}|A|^2\,(R-r)^2 \leq \int_{\cB_R} |\nabla (R-r)|^2 =
    \Area \, (\cB_{R}) \, .
\end{equation}
Consequently, we get an a priori bound for the area of an
intrinsic ball in a stable minimal disk
\begin{equation}
 \Area \, (\cB_{R}) \leq 4\pi \, R^2 / 3 \, .
\end{equation}
(This area bound is the starting point in \cite{CM2}.) We will use
two generalizations of this argument below.  The first will get a
polynomial area bound for  embedded minimal disks with bounded
curvature.  The second generalization will be to bound the area of
a $1/2$-stable sector (a sector is a specific type of subdomain of
an intrinsic ball).

\vskip2mm \noindent {\underline{A curvature estimate for embedded
disks with bounded total curvature}}.   The following curvature
estimate for embedded minimal disks $\Sigma$ generalizes a result
of Schoen and Simon (theorem $2.5$ in \cite{CM1}):

\begin{Pro} \label{p:scsi}
(Corollary $1.18$ in \cite{CM4}). Given $C_I$, there exists $C_P$
so that if
\begin{equation}
    \int_{\cB_{2s}} |A|^2 \leq C_I \, ,
\end{equation}
 then
\begin{equation}
    \sup_{\cB_{s}} |A|^2 \leq C_P \, s^{-2} \, .
\end{equation}
\end{Pro}

 By \eqr{e:jacga}, bounds on $\Area (\cB_{t})/t^2$,
${\text{Length}}(\partial \cB_{t})/t$, or $\int_{\cB_t} |A|^2$ are
equivalent (if we are willing to go to subballs). Therefore,
Proposition \ref{p:scsi} gives an a priori curvature estimate when
any one of these three quantities is bounded. It  is important
that $\Sigma$ is an embedded
 disk; e.g., the catenoid is
complete, has finite total curvature, and is not flat.

\begin{figure}[htbp]
    \setlength{\captionindent}{20pt}
    \begin{minipage}[t]{0.5\textwidth}
\centering\input{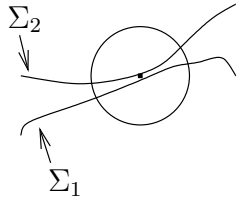}
    \caption{Two sufficiently close disjoint
minimal surfaces with bounded curvatures must each be nearly stable.}
\label{f:7}
\end{minipage}
\end{figure}

\vskip2mm \noindent {\underline{Nearby disjoint surfaces with
$|A|^2 \leq 4$ are nearly stable}}.
 See fig. \ref{f:7}.
We will see that if two disjoint minimal surfaces with bounded
curvature come close enough together, then each of them is
$1/2$-stable.  This $1/2$-stability property is a weakening of
stability which is still sufficiently strong to imply estimates on
the area and curvature of the surfaces.

Before making this precise,
 recall that  the linearization of the minimal graph equation is
  the Jacobi equation (for
 minimal surfaces)
 \begin{equation}
    L u= \Delta u + |A|^2 u = 0 \, .
 \end{equation}
  A domain $\Omega \subset \Sigma$ is
said to be stable if $\int \phi \, L \phi \leq 0$ for every $\phi$
with compact support in $\Omega$; i.e., if we have the stability
inequality
\begin{equation}        \label{e:stabin0}
    \int |A|^2 \, \phi^2 \leq \int |\nabla
\phi|^2 \, .
\end{equation}
We will say that $\Omega$ is $1/2$-stable if we have the weaker
$1/2$-stability inequality
\begin{equation}   \label{e:stabin}
      \int |A|^2 \, \phi^2 \leq 2 \, \int |\nabla \phi|^2 \, .
\end{equation}
One useful criterion for stability is that if $u > 0$ and $Lu = 0$
on $\Omega$, then $\Omega$ is stable; similarly, when $u>0$ and
$Lu/u$ is small, then $\Omega$ is $1/2$-stable (cf. section 2 of
\cite{CM4}).

Using the above criteria for stability, we can now explain why
nearby, but disjoint, minimal disks with bounded curvature must be
$1/2$-stable.   Namely,
  if two disjoint minimal disks with $|A|^2 \leq 4$ come
close at a point, then it is not hard to see that one can be
written as a (normal exponential) graph over the other of a
function $u>0$ with $Lu/u \approx 0$ and, consequently, each is
$1/2$-stable. This is similar to the case of geodesics.
 See lemmas $2.6$  and $2.11$ of \cite{CM4} for the precise statements.

\vskip2mm  Having these tools at our disposal, we turn next to the
main steps in the proof of the existence of the small multi-valued
graphs near a point of large curvature.

\vskip2mm \noindent {\underline{Polynomial area bounds when $|A|^2
\leq 4$}}. In general, the volume comparison theorem from geometry
implies that a surface with bounded curvature has at most
exponential area growth (as is the case for hyperbolic space).
However, we will see that an embedded minimal disk in $\RR^3$ with
bounded curvature actually has polynomial area growth.  This will
be used to prove a doubling property, that is, to find arbitrarily
large balls where the area of the double ball increases by at most
a bounded factor (see Corollary \ref{c:choosing} below for the
precise statement).  Notice that this kind of doubling occurs for
polynomial growth but not for exponential growth, i.e.,
$\lim_{r\to \infty} \frac{(2r)^p}{r^p} = 2^p < \infty$ while
$\lim_{r\to \infty} \frac{\e^{2r}}{\e^r} = \infty$.

 To get this polynomial area bound, we first show that ``most of''
 an embedded minimal disk $\Sigma$ with $|A|^2 \leq 4$ is nearly stable.
 More precisely, the next lemma decomposes $\Sigma$ into a union
 of disjoint $1/2$-stable domains $\Omega_j$ and a remainder with
 bounded area.  Moreover, the lemma also gives a cutoff function
with bounded energy which will be used in the $1/2$-stability
inequality.

\begin{Lem} \label{l:deltstb}
(Lemma $2.15$ in \cite{CM4}).
There exists $C_1$ so:
If $0 \in \Sigma \subset B_{2R}$,
$\partial \Sigma \subset \partial B_{2R}$, and $|A|^2 \leq 4$, then
there exist disjoint
$1/2$-stable subdomains $\Omega_j \subset \Sigma$ and
a function $\chi \leq 1$ which vanishes on $B_R \cap \Sigma \setminus
\cup_j \Omega_j$ so that
\begin{align} \label{e:abd}
\Area ( \{ x \in B_R \cap \Sigma \, | \, \chi (x) < 1 \} ) &
\leq C_1 \, R^3 \, , \\
\int_{\cB_R} |\nabla \chi|^2 & \leq C_1 \, R^3 \, . \label{e:chi}
\end{align}
\end{Lem}

\begin{proof}
(Sketch).
Fix $\rho > 0$ small.
Given $x \in  \Sigma$,
let $\Sigma_{x}$ be the component of $B_{\rho}(x)\cap \Sigma$
with $x \in \Sigma_{x}$ and let $B_x^{+}$ be the component
of $B_{\rho}(x) \setminus \Sigma_{x}$
which $\nn(x)$ points into.
Set, see fig. \ref{f:vb},
\begin{equation} \label{e:bplus}
\vb = \{ x \in B_R \cap \Sigma \, | \,
B_x^{+} \cap \Sigma \setminus
\cB_{4 \, \rho }(x) = \emptyset \}
\end{equation}
and let $\{ \Omega_j \}$ be the components of $B_{R} \cap \Sigma
\setminus \overline{\vb}$. It follows from the previous subsection
that each $\Omega_j$ comes equipped with a positive solution $u_j$
of the minimal graph equation and  is therefore
 $1/2$-stable for $\rho$ sufficiently small
 (it is obvious that the function $u_j$ exists locally, but it requires a slight argument
 to see that it is globally well-defined).

 The function $\chi$ is a linear cutoff function on the
$\rho$-tubular neighborhood of $\vb$, i.e., set
\begin{equation}
\chi (x)=
\begin{cases}
0 & \hbox{ if } x \in \vb \, , \\
\dist_{\Sigma} (x, \vb ) / \rho & \hbox{ if } 0 \leq \dist (x ,
\vb) \leq
    \rho \, , \\
1 & \hbox{ otherwise }\, .
\end{cases}
\end{equation}
Finally, we use a simple ball counting argument to get \eqr{e:abd}
and \eqr{e:chi}.  Roughly speaking, the definition of $\vb$ allows
us to cover it by a collection of extrinsic balls $\{
B_{2\rho}(x_i) \}$ so that the ``half half-balls'' $B_{x_i}^{+}$
are essentially disjoint.  Since each of these disjoint ``half
half-balls'' has volume $\approx \rho^3$ and is contained in a
fixed ball in $\RR^3$, there are at most $C \, (R/\rho)^3$ such
balls. This and the curvature bound easily give \eqr{e:abd} and
\eqr{e:chi}.
\end{proof}

\begin{figure}[htbp]
    \setlength{\captionindent}{20pt}
\begin{minipage}[t]{0.5\textwidth}
\centering\input{blow8.pstex_t}
    \caption{The set $VB$ in \eqr{e:bplus}:
$x \in VB$ and $y \in \Sigma\setminus VB$.}
\label{f:vb}
\end{minipage}
\end{figure}

Using the decomposition from Lemma \ref{l:deltstb} in
\eqr{e:jacga} (as in \eqr{e:jacgab}) gives polynomial area bounds
for intrinsic balls with bounded curvature (Lemma $3.1$ in
\cite{CM4}):

\begin{Pro}     \label{p:tcgr}
 If $0 \in \Sigma \subset B_{2\, R}$, $\partial \Sigma
\subset
\partial B_{2\, R}$, and $|A|^2 \leq  4$, then
\begin{equation} \label{e:tcgr}
 \int_{0}^R
\int_{0}^t \int_{\cB_s} |A|^2 \, ds   \, dt =
2( \Area (\cB_{R})- \pi \, R^2)
\leq 6 \,\pi \, R^2 + 20 \, C_1 \, R^5 \,  .
\end{equation}
\end{Pro}

\begin{proof}
 Let the constant $C_1$, the function $\chi$, and the subdomain $\cup_j \Omega_j$ be
given by Lemma \ref{l:deltstb}.  In particular, the function $\chi
(R-r)$ vanishes off of $\cup_j \Omega_j$. Using $\chi (R-r)$ in
the $1/2$-stability inequality (i.e., in \eqr{e:stabin}), the
absorbing inequality and \eqr{e:chi} give
\begin{align} \label{e:stab1}
\int |A|^2 \chi^2 (R-r)^2 & \leq
2 \, \int \, \left( \chi^2 + 2 \chi \, R \, |
\nabla \chi | + R^2 \, |\nabla \chi|^2 \right)
 \\ & \leq
6 \, R^2 \, \int_{\cB_R} |\nabla \chi|^2
+ 3 \int \chi^2
\leq 6 \, C_1 \, R^5  + 3\, \Area (\cB_R)\, . \notag
\end{align}
On the other hand, combining the area bound \eqr{e:abd} for $\{
\chi < 1 \}$ and $|A|^2 \leq 4$ gives
\begin{equation} \label{e:stab1a}
    \int |A|^2 (1-\chi^2) (R-r)^2  \leq 4 \, R^2 \,  \Area ( \{ x \in
    B_R \cap \Sigma \, | \, \chi (x) < 1 \} )  \leq 4 \, C_1 \, R^5 \,
    .
\end{equation}
We see \eqr{e:tcgr} by using \eqr{e:stab1} and \eqr{e:stab1a} in
\eqr{e:jacga} to get
\begin{equation} \label{e:stab3}
4\,(\Area \, (\cB_{R}) - \pi \, R^2) =
\int |A|^2 (R-r)^2
\leq 10 \, C_1 \, R^5 + 3 \, \Area \, (\cB_R) \, .
\end{equation}
\end{proof}

  Using the polynomial area growth proven in Proposition \ref{p:tcgr},
  it is now standard to find
large intrinsic balls with a fixed doubling for area (and hence
also for total curvature by \eqr{e:jacga}):

\begin{Cor} \label{c:choosing}
(Corollary $3.5$ in \cite{CM4}). There exists $C_2$ so that given
$\beta , R_0 > 1$, we get $R$ so: If $0 \in \Sigma \subset B_{R}$
is an embedded minimal disk, $\partial \Sigma \subset \partial
B_{R}$, $|A|^2(0) = 1$, and $|A|^2 \leq 4$, then there exists $R_0
\leq s < R / (2 \, \beta)$ with
\begin{align}
\int_{\cB_{3 \, s}} |A|^2 &\leq C_{2} \, s^{-2} \, \Area \, (
\cB_{s} ) \,  , \label{e:fe3a} \\
\beta^{-10} \, \int_{\cB_{2 \, \beta \, s}} |A|^2 &\leq C_{2} \,
s^{-2} \, \Area \, ( \cB_{s} ) \, . \label{e:fe3b}
\end{align}
\end{Cor}

\begin{proof} (Sketch)  Since the argument is similar, we sketch
only the proof of \eqr{e:fe3a}.  By \eqr{e:jacga}, it is easy to
see that
\begin{equation}    \label{e:jacgaqq}
\int_{\cB_{3 \, s}} |A|^2 \leq  8 \, s^{-2} \, \Area (\cB_{4 \,
s}) \, .
\end{equation}
Therefore, to prove \eqr{e:fe3a},  it suffices to find $s \geq
R_0$ with
\begin{equation}    \label{e:jacgaq}
\Area (\cB_{4 \, s}) \leq C_2 / 8  \, \Area (\cB_{ s}) \, .
\end{equation}
To do this, we use the bounds for the area given by Proposition
\ref{p:tcgr} to get
\begin{equation} \label{e:fe1}
\left( \min_{1 \leq n \leq m-1} \frac{\Area ( \cB_{4^{n} \, R_0 }
)} {\Area ( \cB_{4^{n-1} \, R_0 } )} \right)^m \leq \frac{\Area (
\cB_{4^{m} \, R_0 } )} {\Area ( \cB_{R_0 } )} \leq C \, 4^{5m} \,
R_0^3 \, .
\end{equation}
Choosing $m$ large so that $C \, R_0^3 \leq 2^m$ and taking the
$m$-th root of \eqr{e:fe1} gives
\begin{equation} \label{e:fe2}
 \min_{1 \leq n \leq m-1} \frac{\Area ( \cB_{4^{n} \, R_0 }
)} {\Area ( \cB_{4^{n-1} \, R_0 } )} \leq 2 \, (4^{5} ) \, .
\end{equation}
Let $j$ be where the minimum in \eqr{e:fe2} is achieved, so that
 we get \eqr{e:jacgaq} with $C_2 = 16 \, ( 4^5)$ and
$s= 4^{j-1} \, R_0$.
\end{proof}

We will not go into why we need both \eqr{e:fe3a} and \eqr{e:fe3b}
since this is somewhat technical.  One consequence of this
doubling property given by \eqr{e:fe3a} is an average curvature
bound for the intrinsic ball $\cB_{3s}$. Together with the
curvature estimate in Proposition \ref{p:scsi}, this implies that
at least ``most of'' $\cB_{3s}$ has a point-wise scale-invariant
curvature bound of the form $|A|^2 \leq C \, r^{-2}$. We will see
later that much more is true; namely, most of $\cB_{3s}$ will be
almost stable (and therefore satisfy even better estimates). To
make this precise, we next discuss a useful way to subdivide
intrinsic balls into intrinsic sectors.

\vskip2mm \noindent {\underline{Area and total curvature of
$1/2$-stable sectors}}.  We next define intrinsic sectors and give
estimates for the area and total curvature of a $1/2$-stable
intrinsic sector.  Given a curve $\gamma \subset  \partial
\cB_{s}$, define the intrinsic (truncated) sector, see fig.
\ref{f:9},
\begin{equation} \label{e:defsip}
  S = S_{R}(\gamma) = \{ \exp_0 (v) \, | \,  s \leq |v| \leq
s + R {\text{ and }}
\exp_0 (s \, v / |v|) \in   \gamma \} \, .
\end{equation}
(Here $\exp_0: \RR^2 \to \Sigma$ is the exponential map at $0$.)
The simplest example of a sector is when $\Sigma = \RR^2$ and
$\gamma$ is given in polar coordinates as $\{ \rho = s , \,
\theta_1 \leq \theta \leq \theta_2 \}$; in this case, the sector
$S_{R}(\gamma)$ is just $\{ s \leq \rho \leq s + R , \, \theta_1
\leq \theta \leq \theta_2 \}$.  We will often refer to $\gamma$ as
the inner boundary of the sector and the geodesic rays through
$\partial \gamma$ as the sides of the sector.

Just as stability led to area bounds for intrinsic balls in
\eqr{e:jacgab}, we can bound the area of a $1/2$-stable sector.
However, to make the cutoff function compactly supported on a
sector, we must also cutoff along the inner boundary  and along
the sides.  This introduces new terms in the upper bound for the
area. The next lemma, which is an easy consequence of lemma II.1.1
 and remark II.1.32 in \cite{CM3},  gives the resulting estimate
 ($k_g$ is the geodesic
curvature of the curve $\gamma$).

\begin{Lem} \label{l:stable}
\cite{CM3}.
Suppose that $S$ is $1/2$-stable,
$\int_{\gamma} (1 + |k_g| \, s) \leq C_0 \, m \, s$,
 $R > 2 \, s$, and  for $x \in S$ we have
$\sup_{\cB_{s/4}(x)}|A|^2 \leq C_0  s^{-2}$.
Then for  $\Omega > 2$ and $2 \, s \leq t \leq 3 R / 4$
\begin{align}
t \, \int_{\gamma} k_g \leq {\text{Length}} (\partial \cB_t \cap
S) &\leq C_3 \, (m +R/s) \, t \, , \label{e:lbound} \\
\int_{ S_{R/\Omega} \setminus \cB_{\Omega s} }
 |A|^2 &\leq C_1 \, R
/ s + C_2 \, m / \log \Omega \, .  \label{e:cbound}
\end{align}
\end{Lem}

\begin{proof}
(Sketch). First, as in \eqr{e:jacgab}, we combine the
$1/2$-stability inequality and the analog of \eqr{e:jacga} for
intrinsic sectors to get \eqr{e:lbound}. The $R/s$ term in
\eqr{e:lbound} comes from cutting off linearly on a (roughly)
$s$-tubular neighborhood of the sides of the sector (which have
length $R$). Once we have the quadratic area growth given by
\eqr{e:lbound}, we can then use a radial logarithmic cutoff in the
stability inequality to get \eqr{e:cbound}.
\end{proof}

When we apply Lemma \ref{l:stable}, the bound $\sup_{ \cB_{s/4}(x)
}|A|^2 \leq C_0 s^{-2}$ on $S$  will be given  by starting with a
slightly larger $1/2$-stable sector and then coming in from the
boundary (using a comparison theorem to guarantee that the
required intrinsic balls are contained in the larger sector).

Roughly speaking, Lemma \ref{l:stable} shows  that if we can find
a large $1/2$-stable sector, then we can find subsectors with
small average curvature.  We next explain how this leads to
finding a subsector with small total curvature.

\vskip2mm \noindent {\underline{Small total curvature of stable
subsectors}}. Given a  $1/2$-stable sector $S$ as in Lemma
\ref{l:stable} over a long curve $\gamma$, we can subdivide
$\gamma$ into subcurves so that one of the sectors over these
subcurves has small total curvature.    To see this, suppose that
${\text{Length}}(\gamma) \approx m \, s$ for some $m$ large. If
$m/\log \Omega$ is larger than $R/s$, then \eqr{e:cbound} says
that the truncated sector has small average curvature
\begin{equation}    \label{e:smac}
    \int_{S_{R/\Omega} \setminus \cB_{\Omega s}}|A|^2 \leq C \, m
/ \log \Omega \, .
\end{equation}
 If we now subdivide $\gamma$ into $m$ subcurves, it
follows from \eqr{e:smac} that the truncated sector $\hat{S}$ over
one of these subcurves has small total curvature
\begin{equation}
    \int_{\hat{S}} |A|^2 \leq C / \log \Omega \, .
\end{equation}

\vskip2mm \noindent {\underline{The $N$-valued graph}}. The next
step is to use the $1/2$-stable sectors with small total curvature
to construct the $N$-valued graph $\Sigma_g$.   To see this,
suppose that $S$ is as in Lemma \ref{l:stable} and $\hat{S}
\subset S$ is a subsector with small total curvature
\begin{equation}
    \int_{\hat{S}} |A|^2 \leq \epsilon \, ,
\end{equation}
 where $\epsilon > 0$ can be chosen small by taking $\Omega$ and $m$ large.
The small total curvature (and stability) implies that
\begin{equation}        \label{e:tnvg}
    |A|^2 \leq C \, \epsilon \, r^{-2} {\text{ on }} \hat{S} \, .
\end{equation}
We will use this estimate \eqr{e:tnvg} to build out $\Sigma_g$.
 For each point
  $x$   on the inner boundary of $\hat{S}$, let
$\gamma_x \subset \Sigma$ be the geodesic leaving $x$ orthogonally
to $\partial \hat{S}$.  Integrating the curvature bound
\eqr{e:tnvg} along $\gamma_x$ implies first that, as a curve in
$\RR^3$, $\gamma_x$ has small total geodesic curvature and is
close to a line segment. Second, it also implies that the unit
normal $\nn$ to $\Sigma$ has small oscillation along $\gamma_x$.
It is now easy to see that as we vary the initial point $x$ on the
inner boundary of $\hat{S}$, the geodesics $\gamma_x$ trace out a
multi-valued graph $\Sigma_g$.

\vskip2mm Finally, we review how these steps fit together to prove
Theorem \ref{t:blowupwindinga}.

\vskip2mm \noindent {\bf{Sketch of the proof of Theorem
\ref{t:blowupwindinga}}}: After rescaling by $C/r_0$,  we can
assume that
\begin{equation}    \label{e:afterscale}
    \sup_{B_{C} \cap \Sigma}|A|^2\leq
    4\,|A|^2(0)= 4 \, .
\end{equation}
 The key for proving  Theorem
\ref{t:blowupwindinga} is to find many large intrinsic sectors
with a quadratic curvature bound. To do this,  we first fix some
large $\bar{m}$ and use Proposition \ref{p:scsi} and the lower
bound $|A|^2(0) = 1$ in \eqr{e:afterscale} to get $R_0$ so that
for $R \geq R_0$
\begin{equation}
    {\text{Length}}(\partial \cB_R) \geq \bar{m} \, R \, .
\end{equation}
Using the upper bound $|A|^2 \leq 4$ in \eqr{e:afterscale}, we can
apply Corollary \ref{c:choosing} to get $R_3
>R_0$
and many long disjoint curves $\tilde{\gamma}_i \subset
\partial \cB_{R_3}$ so the sectors over $\tilde{\gamma}_i$ have
uniformly bounded total curvature $\int |A|^2$.   Proposition
\ref{p:scsi} then gives a (point-wise)
 quadratic curvature  bound on these sectors.

Now that we have these many disjoint sectors with a quadratic
curvature bound (as many as we want by taking $\bar{m}$ large),
two must be close together and, hence, also $1/2$-stable. Lemma
\ref{l:stable} gives a subsector with small total curvature which
then must contain the $N$-valued graph $\Sigma_g$ (this is done in
corollary II.1.34 of \cite{CM3}).

\section{The estimate between the sheets and the
extension of multi--valued graphs}   \label{s:sheets}

In this section, we will give an overview of the proof of
Theorem  \ref{t:spin4ever2} where we extend the multi-valued graphs
almost to the boundary corresponding to B. in fig. \ref{f:f3a}.
This result is significantly more subtle than
the local existence result discussed in the previous section.  We only give a
very rough outline of the proof and refer to \cite{CM3} for the full story.

\begin{figure}[htbp]
    \setlength{\captionindent}{20pt}
    \begin{minipage}[t]{0.5\textwidth}
\centering\input{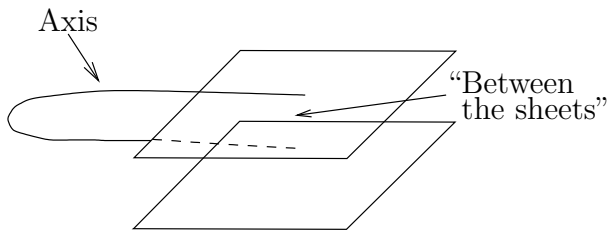}
    \caption{The estimate between the sheets; Theorem \ref{t:tslab}.}
\label{f:f22a}
\end{minipage}
\end{figure}

\begin{figure}[htbp]
    \setlength{\captionindent}{20pt}
    \begin{minipage}[t]{0.5\textwidth}
\centering\input{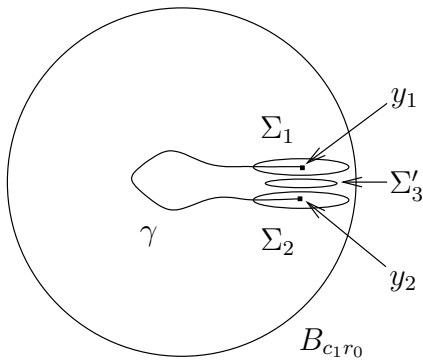}  
    \caption{$y_1$, $y_2$, $\Sigma_1$, $\Sigma_2$, $\Sigma_3'$, and $\gamma$ in
    Theorem \ref{t:tslab}.}
\label{f:f22}
\end{minipage}
\end{figure}

The key component of the proof of the extension theorem is a curvature
estimate ``between the sheets'' for embedded minimal disks in
$\RR^3$.
We will think of an axis (see fig. \ref{f:f22})
for such a disk $\Sigma$ as a point or
curve away from which the
surface locally (in an extrinsic ball) has more than one
component.  With this weak notion of an axis, our estimate is that if
one component of $\Sigma$ is sandwiched between two
others that connect to an axis, then the one that
is sandwiched has curvature estimates; see Theorem \ref{t:tslab}.
The example to keep in
mind is a helicoid and the components are
``consecutive sheets'' away from the axis.

Let  $\gamma_{p,q}$ denote the line segment from $p$ to
$q$ and $T_s(\gamma_{p,q})$ its $s$-tubular neighborhood.
A curve
$\gamma$ is $h$-{\it{almost monotone}} if given $y \in \gamma$,
then $B_{4 \, h}(y) \cap \gamma$ has only one component which
intersects $B_{2\, h }(y)$.
Our curvature estimate ``between the sheets'' is:

\begin{Thm} \label{t:tslab}
(Theorem I.0.8 in \cite{CM3}).   See fig. \ref{f:f22}.
There exist $c_1 \geq 4$, $2 c_2 < c_4 < c_3 \leq 1$ so:
Let $\Sigma \subset B_{c_1 \,r_0}$ be an embedded minimal disk with
$\partial \Sigma \subset
\partial B_{c_1 \, r_0}$ and $y \in \partial B_{2\, r_0}$.
Suppose $\Sigma_1 , \Sigma_2 , \Sigma_3$ are distinct
components of $B_{ r_0}(y) \cap \Sigma$ and $\gamma \subset (
B_{r_0} \cup T_{c_2 \,r_0}(\gamma_{0,y}) ) \cap \Sigma$ is a curve
with $\partial \gamma = \{ y_1 , y_2 \}$ where $y_i \in B_{c_2 \,
r_0 }(y) \cap \Sigma_i$ and each component of $\gamma \setminus
B_{r_0}$ is $c_2 \, r_0$-almost monotone. Then any component
$\Sigma_3'$ of $B_{c_3 \, r_0 }(y) \cap \Sigma_3$ with $y_1 ,
y_2$ in distinct components of $B_{c_4 \, r_0}(y) \setminus
\Sigma'_3$ is a graph.
\end{Thm}

The idea of the proof of Theorem \ref{t:tslab}
 is to show that if this were not the case, then
we could find an embedded stable disk that would be almost
flat and lie in the complement of the original  disk. In fact, we
can choose the stable disk to be sandwiched between the two
components as well. The flatness would force the stable disk to
eventually cross the axis in the original disk, contradicting that
they were disjoint.

The curve $\gamma$ in Theorem
\ref{t:tslab} which plays the role of an axis and connects
$\Sigma_1$ and $\Sigma_2$ can in many instances
be extended using the maximum principle
once it occurs on a given small scale; cf. lemma I.0.11 of \cite{CM3}.
This is used both in the proof of Theorem
\ref{t:tslab} and in the applications of it.

\vskip2mm \noindent {\bf{Overview of the proof of Theorem
\ref{t:spin4ever2}}}: The first step (part I of \cite{CM3}) is to
show Theorem \ref{t:tslab} when the surface is in a slab; i.e.,
when $\Sigma \subset \{ |x_3| \leq \beta \, h \}$.

The second step (part II of \cite{CM3}) is to show that certain
stable disks starting off as multi-valued graphs remain
multi-valued graphs as the outer radius increases. This is needed
when we generalize the results of step one to when the surface is
not anymore in a slab.

Two facts go into the proof of the extension of stable disks.
First, we show that if an almost flat multi-valued graph sits
inside a stable disk, then the outwardly defined intrinsic sector
from a curve which is a multi-valued graph over a circle has a
subsector which is almost flat; cf. Lemma \ref{l:stable} and the
paragraph ``Small total curvature of stable subsectors'' above. As
the  initial multi-valued graph becomes flatter and the number of
sheets in it goes up, the subsector becomes flatter. The second
fact is that the separation between the sheets grows sublinearly;
\eqr{e:slg}.

The first application of these two facts
is to extend the middle sheet as a
multi-valued graph.  This is done by dividing the initial
multi-valued graph (or curve in the graph that is itself a
multi-valued graph over the circle) into three parts where the
middle sheet is the second part.  The idea is then that the first
and third parts have subsectors which are almost flat multi-valued
graphs and the middle part (which has curvature estimates since it
is stable) is sandwiched between the two others.  Hence its sector
is also almost flat.

A thing that adds some technical complications to the above
picture is that in the result about almost flat
subsectors it is important that the ratio between the size of the
initial multi-valued graph and how far one can go out is fixed.
This is because the estimate for the subsector comes from a total
curvature estimate which is in terms of this ratio (see
\eqr{e:cbound}) and can only be made small by looking at a fixed
large number of rotations for the graph. This forces us to
successively extend the multi-valued graph.  The issue is then to
make sure that as we move out in the sector and repeat the
argument we have essentially not lost sheets.  This is taken care
of by using the sublinear growth of the separation between the
sheets together with the Harnack inequality
and the maximum principle. (The maximum
principle is used to make sure that as we try to recover sheets
after we have moved out then we don't hit the boundary of the disk
before we have recovered essentially all of the sheets that we
started with.) The last thing which is used
is theorem $3.36$ in \cite{CM10}.  This is used
to guarantee that, as we patch together these multi-valued
graphs coming from different scales, then the surface that we get
is still a multi-valued graph over a fixed plane.

The third step (part III of \cite{CM3}) is to generalize the curvature
estimate between the sheets to the case where the surface is not anymore
in a slab.  This uses the extension of the stable graphs from step two.

Finally, using steps one,
two, and three we showed Theorem \ref{t:spin4ever2}
in part IV of \cite{CM3}.

\part{The proof of the one-sided curvature estimate}

In appendix A of \cite{CM6} we showed curvature estimates for minimal
hyper-surfaces in $\RR^n$ which on all
sufficiently small scales lie on one side of, but come close to,
a hyper-plane. Such a scale invariant version of Theorem \ref{t:t2} can
(unlike Theorem \ref{t:t2}) be proven quite easily by a blow up argument
using the minimal surface equation.
Moreover, such a scale invariant condition is very similar
to the classical Reifenberg property.  (After all, a subset of $\RR^n$
has the Reifenberg property if it is close on all scales to a hyper-plane.)
As explained in Part \ref{p:p1} (in particular Corollary \ref{c:conecor}),
the significance of Theorem \ref{t:t2}
is indeed that it only requires closeness on one scale.  On the other
hand, this is what makes it difficult to prove and requires us to use
the results discussed in the previous parts together with results
from \cite{CM5}.

Let us briefly outline the proof of the one-sided; i.e.,
Theorem \ref{t:t2}. Suppose that $\Sigma$
is an embedded minimal disk in the half-space $\{x_3> 0\}$. We prove the
curvature estimate by contradiction; so suppose that $\Sigma$ has low
points with large curvature.
Starting at such a point, we decompose
$\Sigma$ into disjoint multi-valued graphs using the existence
of nearby points with large curvature given by corollary III.3.5 of \cite{CM5},
 a blow up argument,
and \cite{CM3}, \cite{CM4}. The key point is then to show
(see Proposition \ref{p:lift} below)
that we can in fact find such a decomposition where the ``next''
multi-valued graph starts off a definite amount below where the
previous multi-valued graph started off. In fact, what we show
is that this definite amount is a fixed fraction of the distance
between where the two graphs started off. Iterating this
eventually forces $\Sigma$ to have points where $x_3<0$.
This is the desired contradiction.

\begin{figure}[htbp]
    \setlength{\captionindent}{20pt}
    \begin{minipage}[t]{0.5\textwidth}
    \centering\input{unot10.pstex_t}
    \caption{Proposition \ref{p:lift}:
Two consecutive blow up points satisfying \eqr{e:pairs}.}
\label{f:f23}
    \end{minipage}
\end{figure}

\begin{Pro}   \label{p:lift}
(See proposition III.2.2 in \cite{CM6} for the precise statement).
See fig. \ref{f:f23}.
There exists $\delta>0$ such that if $(0,s)$
satisfies \eqr{e:pairs} and $\Sigma_0\subset \Sigma$ is the corresponding
(to $(0,s)$) $2$-valued graph over $D_R\setminus D_s$, then
we get $(y,t)$ satisfying \eqr{e:pairs} with $y\in \cone_{\delta}(0)
\cap \Sigma \setminus B_{Cs/2}$ and where $y$ is below $\Sigma_0$.
\end{Pro}

\begin{figure}[htbp]
    \setlength{\captionindent}{20pt}
    \begin{minipage}[t]{0.5\textwidth}
    \centering\input{unot11.pstex_t}
    \caption{Between two consecutive blow up points satisfying  \eqr{e:pairs}
there are a bunch of blow up points satisfying
 \eqr{e:defc1ii}.}  \label{f:f24}
    \end{minipage}\begin{minipage}[t]{0.5\textwidth}
    \centering\input{unot12.pstex_t}
    \caption{Measuring height.  Blow up points and
corresponding multi-valued graphs.}  \label{f:f25}
    \end{minipage}
\end{figure}

To prove this key proposition
(Proposition \ref{p:lift}), we use two decompositions
and two kinds of blow up points.
The first decomposition uses
the more standard blow up points given as pairs $(y,s)$
where $y\in \Sigma$ and $s>0$ is such that
\begin{equation}   \label{e:defc1ii}
\sup_{\cB_{8s}(y)}|A|^2\leq 4|A|^2(y)=4\,C_1^2\,s^{-2}\, .
\end{equation}
The point about such a pair $(y,s)$ is that by \cite{CM3}, \cite{CM4}
(and an argument in
section II.2 of \cite{CM6}
which allows us replace extrinsic balls by
intrinsic ones), then
$\Sigma$ contains a multi-valued graph near $y$ starting
off on the scale $s$.  (This is assuming that $C_1$ is a sufficiently
large constant given by \cite{CM3}, \cite{CM4}.)  The second kind of blow up
points are the ones where (except for a technical issue)
$8$ is replaced by some really large
constant $C$, i.e.,
\begin{equation}  \label{e:pairs}
\sup_{\cB_{Cs}(y)}|A|^2\leq 4|A|^2(y)=4\,C_1^2\,s^{-2}\, .
\end{equation}
The point will then be that we can find blow up points
satisfying \eqr{e:pairs} so that the distance between them is proportional
to the sum of the scales.  Moreover, between consecutive
blow up points satisfying \eqr{e:pairs}, we
can find a bunch of blow up points satisfying \eqr{e:defc1ii};
see fig. \ref{f:f24}.  The advantage
is now that if we look between blow up points satisfying \eqr{e:pairs}, then
the height of the multi-valued graph given by such a pair grows
like a small power of the distance whereas the separation between the sheets
in a multi-valued graph given by \eqr{e:defc1ii} decays like a small
power of the
distance; see fig. \ref{f:f25}.
Now thanks to that the number of blow up points satisfying
\eqr{e:defc1ii} (between two consecutive blow up points satisfying
\eqr{e:pairs}) grows
almost linearly then, even though the height of the graph coming from
the blow up point satisfying \eqr{e:pairs} could move up (and thus work
against us), then the sum of the separations of the graphs coming from the
points satisfying \eqr{e:defc1ii} dominates the other term.  Thus the
next blow
up point satisfying \eqr{e:pairs} (which lies below all the other graphs)
is forced to be a definite amount lower than the previous blow up
point satisfying \eqr{e:pairs}.

\part{The local case - when singular limit laminations can occur}
\label{p:local}

In this part we discuss the differences  between the so-called
local and global cases.  The local case is where we have a
sequence of embedded minimal disks in a ball of fixed radius in
$\RR^3$ - the global case is where the disks are in a sequence of
expanding balls with radii tending to infinity. The main
difference between these cases is that in the local case we can
get limits with singularities. In the global case this does not
happen because in fact any limit is a foliation by flat parallel
planes (assuming that the curvatures of the sequence are blowing
up).

To precisely define the local and global cases, suppose
$\Sigma_i\subset B_{R_i}=B_{R_i}(0)\subset \RR^3$ with $\partial
\Sigma_i\subset
\partial B_{R_i}$ is a sequence of (compact) embedded minimal
disks and
either:\\
(a) $R_i$ equal to a finite constant.\\
(b) $R_i\to\infty$.\\
Case (a) is what we call the {\it local case} and (b) is what we
refer to as the {\it global case}; Theorem \ref{t:t0.1} dealt with
the global case. Recall
that a surface $\Sigma \subset \RR^3$ is
said to be properly embedded if it is embedded and the
intersection of $\Sigma$ with any compact subset of $\RR^3$ is
compact.  We say that a lamination or foliation is proper if each
leaf is proper.

\begin{figure}[htbp]
    \setlength{\captionindent}{4pt}
    \begin{minipage}[t]{0.5\textwidth}
\centering\input{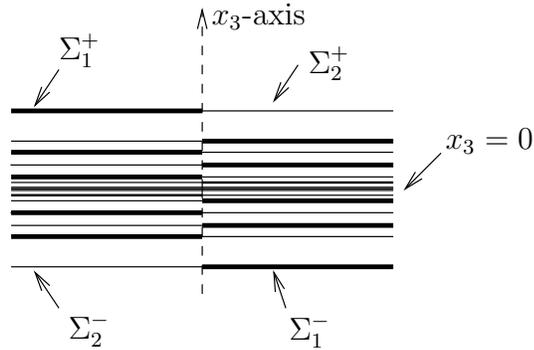}
    \caption{A schematic picture of the limit in Theorem \ref{t:mainb}
    which is not smooth and
     not proper (the dotted $x_3$-axis is part of the limit).
The limit contains four multi-valued graphs
joined at the $x_3$-axis;  $\Sigma_1^+$, $\Sigma_2^+$ above
the plane $x_3=0$ and $\Sigma_1^-$, $\Sigma_2^-$ below the plane.  Each of
the four spirals into the plane.
}\label{f:f26}
\end{minipage}
\end{figure}

To explain the difference between the two cases,  we consider
sequences of minimal disks $\Sigma_i$ as above where the
curvatures blow up, e.g.,
\begin{equation}
    \lim_{i\to \infty} \sup_{B_1 \cap \Sigma_i} |A|^2 = \infty \,
    .
\end{equation}
(Of course, if the curvatures do not blow up, then the
Arzela-Ascoli theorem easily gives smooth convergence of a
subsequence in either case.)
 In the global case, Theorem \ref{t:t0.1} gives a
subsequence of the $\Sigma_i$ converging off of a Lipschitz curve
to a foliation by parallel planes.  In particular, the limit is a
(smooth) foliation which is proper.  However,
 we showed in \cite{CM16} (see Theorem \ref{t:mainb} below) that
smoothness and properness of the limit can fail in the local case;
 cf. fig. \ref{f:f26}.

To illustrate the key issue for the failure of properness, suppose
that $|A|^2 (0) \to \infty$ as $i\to \infty$. In either the local
or global case (see Section \ref{s:13}), we get a sequence of
$2$-valued graphs which converges to a minimal graph $\Sigma_0$
through $0$ (this graph is a plane in the global case).
Furthermore, by the one-sided curvature estimate (see Corollary
\ref{c:conecor}), the intersection of $\Sigma_i$ with a low cone
about $\Sigma_0$ consists of multi-valued graphs for $i$ large.
There are now two possibilities:
\begin{itemize}
\item
The multi-valued graphs in this low cone close up in the limit.
\item
The limits of these multi-valued graphs spiral infinitely into
$\Sigma_0$.
\end{itemize}
In the first case, where properness holds, the sequence converges
to a foliation in a neighborhood of $0$.  In the second case,
where properness fails, the sequence converges to a lamination
away from $0$ but cannot be extended smoothly to any neighborhood
of $0$.   The proof of properness in the global case   is
described in Section \ref{s:proper}.

 In the local case, Theorem \ref{t:mainb} constructs  a sequence of disks
$\Sigma_i \subset B_1$ as above where
 the curvatures blow up only at $0$ and
$\Sigma_i \setminus \{ {\text{$x_3$-axis}} \}$ consists of two
multi-valued graphs for each $i$; see (1), (2), and (3).
Furthermore (see (4)),
$\Sigma_i  \setminus \{ x_3 = 0 \}$  converges to
two
  embedded minimal disks $\Sigma^- \subset \{ x_3 < 0 \}$
and  $\Sigma^+ \subset \{ x_3 > 0 \}$ each of which spirals
 into $\{ x_3 =
0 \}$ and thus is not proper; see fig. \ref{f:f26}.

\begin{Thm} \label{t:mainb}
(Theorem 1 in \cite{CM16}).  See fig. \ref{f:f26}.
There is  a sequence  of compact embedded minimal disks $0 \in
\Sigma_i \subset  B_1 \subset \RR^3$
 with $\partial \Sigma_i \subset
\partial B_1$ and
containing the vertical segment
$\{ (0,0,t) \, | \, |t|<1 \} \subset \Sigma_i$ so:\\
(1) $\lim_{i\to \infty} |A_{\Sigma_i}|^2 (0) = \infty$.\\
(2) $\sup_i \sup_{\Sigma_i \setminus B_{\delta} } |A_{\Sigma_i}|^2 < \infty$
for all $\delta > 0$.\\
(3)  $\Sigma_i \setminus \{ {\text{$x_3$-axis}} \} =
\Sigma_{1,i} \cup \Sigma_{2,i}$ for multi-valued graphs
$\Sigma_{1,i}$ and $\Sigma_{2,i}$\\
(4) $\Sigma_i  \setminus \{ x_3 = 0 \}$  converges to     two
embedded minimal disks $\Sigma^{\pm} \subset \{ \pm x_3 > 0 \}$
with $\overline{\Sigma^{\pm}} \setminus \Sigma^{\pm} = B_1 \cap \{
x_3 = 0\}$. Moreover, $\Sigma^{\pm} \setminus \{
{\text{$x_3$-axis}} \} = \Sigma_1^{\pm} \cup \Sigma_2^{\pm}$ for
multi-valued graphs $\Sigma_{1}^{\pm}$ and $\Sigma_{2}^{\pm}$
which spiral into $\{ x_3 =
0 \}$; see fig. \ref{f:f26}.
\end{Thm}

It follows from
(4) that $\Sigma_i \setminus \{ 0 \}$ converges to a lamination
of $B_1 \setminus \{ 0 \}$    (with leaves $\Sigma^-$, $\Sigma^+$,
and $B_1 \cap \{ x_3 = 0 \} \setminus \{ 0 \}$) which
 does not extend to a lamination of $B_1$.
Namely, $0$ is not a removable singularity.

\vskip2mm  The example in Theorem \ref{t:mainb} was constructed
 using the Weierstrass representation. Recall that if $\Omega
\subset \CC$ is a domain, then the classical Weierstrass
representation starts from a meromorphic function $g$ on $\Omega$
and a holomorphic one-form $\phi$ on $\Omega$ and associates a
(branched) conformal minimal immersion $F: \Omega \to \RR^3$ by
\begin{equation}    \label{e:ws1}
    F(z) = {\text{Re }} \int_{\zeta \in \gamma_{z_0,z}}
\left( \frac{1}{2} \, (g^{-1} (\zeta) - g (\zeta) )
    , \frac{i}{2} \, (g^{-1} (\zeta)
    +g (\zeta) ) , 1 \right) \, \phi (\zeta) \, .
\end{equation}
Here $z_0 \in \Omega$ is a fixed base point and the integration is
along a path $\gamma_{z_0,z}$ from $z_0$ to $z$. The choice of
$z_0$ changes $F$ by adding a constant. When $g$ has no zeros or
poles and $\Omega$ is simply connected, then
 $F(z)$ does not depend on the choice of path
$\gamma_{z_0,z}$.

\vskip2mm \noindent {\bf{Sketch of the proof of Theorem
\ref{t:mainb}}}:   We construct a one-parameter family (with
parameter $a\in (0,1/2)$) of minimal immersions $F_a$ by making a
specific choice of Weierstrass data $g = \e^{i h_a}$ (where $h_a =
u_a+i\, v_a$), $\phi = dz$, and domain $\Omega_a$ to use in
\eqr{e:ws1}. Namely,
  for each $0<a<1/2$,  we define
\begin{equation}   \label{e:dg0}
    h_a (z)  = \frac{1}{a} \,
\arctan \left( \frac{z}{a} \right) {\hbox{ on }}
    \Omega_a = \{ (x,y) \, | \, |x| \leq 1/2 ,
\, |y| \leq (x^2 + a^2)^{3/4}/2  \}
 \, .
\end{equation}
Note that the function $h_a$ is well-defined since $\Omega_a$ is
simply connected and $\pm i \, a \notin \Omega_a$.

It remains to verify that this sequence of minimal immersions has
the desired properties.
 Properties (1) and (2) in Theorem \ref{t:mainb}
 follow easily from calculating the curvature
$K_a$
\begin{equation}
    \K_a (z) = \frac{-|\partial_z h_a|^2}{\cosh^4 v_a}
    =  \frac{-|z^2 + a^2|^{-2}}{\cosh^{4} \left( {\text{Im }}
    \arctan ( z/a) / a \right)}  \, .  \label{e:K2}
\end{equation}
Properties (3) and (4) as well as the convergence away from $0$
follow rather easily by integrating the equations \eqr{e:ws1}. The
key point in the proof, and only remaining point, is to show that
the immersions $F_a : \Omega_a \to \RR^3$ are embeddings.  We do
this by showing that each vertical line segment in the domain
$\Omega_a$ is mapped by $F_a$ to curve in a horizontal plane
(i.e., a plane where $x_3$ is constant) and each such curve is a
graph over a fixed line segment; see \cite{CM16} for the details.

\end{document}